\documentclass[12pt]{article}
\usepackage{graphicx}
\usepackage[outdir=./Figures/]{epstopdf}
\usepackage{amsmath}
\usepackage{amssymb}
\usepackage{multirow}
\usepackage{color}
\renewcommand{\epsilon}{\varepsilon}

\usepackage{amsmath}
\usepackage{amssymb}
\usepackage[ansinew]{inputenc}
\usepackage{graphicx}
\usepackage{hyperref}

\usepackage[authoryear]{natbib}
\newtheorem{theorem}{Theorem}[section]
\newtheorem{corollary}{Corollary}[section]
\newtheorem{lemma}{Lemma}[section]
\newtheorem{example}{Example}[section]
\newtheorem{remark}{Remark}[section]

\def\IM{{\bf I\kern-.25em M}}

\def\3{\ss}

\newcommand{\bea}{\begin{eqnarray*}}
\newcommand{\eea}{\end{eqnarray*}}
\newcommand{\be}{\begin{eqnarray}}
\newcommand{\ee}{\end{eqnarray}}

\newcommand{\ba}{\begin{array}}
\newcommand{\ea}{\end{array}}
\def\3{\ss}

\graphicspath{ {./Figures/} }

\begin{document}

\title{Optimal designs for comparing regression curves -  dependence within and between groups }

\author{
 {\small Kirsten Schorning}
\\
 {\small Technische Universit\"at Dortmund}
 \\
 {\small Fakult\"at Statistik}
 \\
 {\small 44221 Dortmund, Germany}
 \and
  {\small Holger Dette} \\
{\small Ruhr-Universit\"at Bochum} \\
{\small Fakult\"at f\"ur Mathematik} \\
{\small 44780 Bochum, Germany}
}
\date{}
\maketitle

\date{}

\maketitle

 \begin{abstract}
We consider the problem of designing experiments for the comparison of two regression curves describing the relation between
a predictor and a response in two groups, where the data
between and within the group may be dependent. In order to derive efficient designs we use results from    stochastic analysis to identify the best linear unbiased estimator (BLUE) in a corresponding continuous time model.
It is demonstrated that in general simultaneous estimation
using the data from  both groups  yields more precise results
than estimation of the parameters separately  in the two groups. Using the BLUE from simultaneous estimation, we then construct an efficient linear estimator   for finite sample size by minimizing the mean squared error between the optimal solution in the continuous time model and its discrete approximation with respect to the weights (of the linear estimator).  Finally, the optimal design points are  determined by minimizing the maximal width of a simultaneous  confidence band for the difference of the two regression functions.
 The advantages of the new approach  are illustrated by means of a simulation study, where it is shown that the use of the optimal designs yields substantially narrower confidence bands than the application of uniform designs.

 \end{abstract}



 \section{Introduction}
\label{sec1}
\def\theequation{1.\arabic{equation}}
\setcounter{equation}{0}

The application of optimal or efficient designs can improve the accuracy of statistical analysis substantially and meanwhile
there exists a well established  and powerful  theory for the construction of (approximate) optimal designs for independent observations, see
for example the monographs of  \cite{pukelsheim2006} or \cite{fedleo2013}.
In contrast, the determination of  optimal  designs for efficient statistical analysis  from dependent data  is more challenging because
the corresponding optimization problems are in general not convex and therefore the powerful tools of convex analysis are not
applicable.  Although design problems for correlated data have been discussed for a long time
   \citep[see, for example][who use asymptotic arguments
 to develop  continuous but in general   non-convex optimzation problems in this context]{sackylv1966,sackylv1968,bickherz1979,N1985a}
  a  large part  of the discussion is restricted
to  models  with a small number of parameters and we refer
  \cite{pazmue2001},  \cite{muepaz2003}, \cite{detkunpep2008}, \cite{KisStehlik2008}, \cite{harstu2010}, \cite{RODRIGUEZDIAZ2017287}, \cite{campos2015} and
  \cite{attia2020}
among others.

  Recently,  \cite{DetPZ2012}  suggest a more systematic approach
  to the problem and
   determine (asymptotic) optimal designs for least squares estimation,  under the additional assumption
 that the regression functions are eigenfunctions of an integral operator associated with the covariance kernel of the error process.
This approach  refers to  models, where the  regression functions  are  eigenfunctions of the
integral operator corresponding to the
 covariance kernel, which is used to describe the dependencies.
For more general models
\cite{detpepzhi2015}  propose
  to   construct  the  optimal design and estimator simultaneously. More precisely,
  they construct a class of estimators and corresponding optimal designs with a variance converging (as the sample size increases) to the optimal variance in the continuous time model.
   \cite{detkonzhi2017} propose  an alternative strategy  for this purpose. They start with the  construction   of the best linear unbiased estimator (BLUE) in the continuous time model
using stochastic calculus and  determine in a second step  an implementable design, which is ``close'' to the solution in the continuous time model.
 By this approach these authors  are able to  provide an easily implementable  estimator with a corresponding design which is  practically non distinguishable from the weighted least squares estimate (WLSE)  with corresponding optimal design.
 Their results are applicable for a  broad class of linear regression models with various covariance kernels and have recently been extended to the situation, where
 also derivatives of the process can be observed \citep[see ][]{dette2019}.

\cite{dette2016} and   \cite{DETTE2017273}  propose designs for the comparison of  regression curves from two independent samples, where the latter reference
also allows for dependencies  within the samples.  Their work is motivated  by applications in drug development, where a
 comparison between  two regression models that describe the relation between a common response and the same covariates for two groups
is  used to establish the non-superiority of one model to the other or to check whether the difference between two regression models  can be neglected.
For example, if the  similarity between two regression functions describing the
dose response relationships in the groups individually has been established
subsequent inference  in drug development  could be based on the combined samples such that a more efficient statistical
analysis is possible on the basis of the larger population.  Because of its importance several procedures for the comparison of curves
have been investigated in linear and nonlinear models \citep[see][among others]{liubrehaywynn2009,gsteiger2011,liujamzhang2011,detmolvolbre2015,bretzmoelldette2016,moellenhoff2018,mollenhoff2019}.
 Designs minimizing the maximal
 width of a (simultaneous) confidence  band for the difference between the regression curves calculated from two independent groups 
 are determined by \cite{dette2016} and   \cite{DETTE2017273}, who also demonstrate that
 the use of these designs yields to substantially narrower confidence bands.

 While these  results  refer to independent groups it is the purpose of the present paper to investigate designs  for the comparison of  regression curves corresponding to two groups, where
 the data within the groups and  between the groups may be dependent. It will be demonstrated
 that  in most cases simultaneous estimation of the parameters in the regression models using the data from both groups  yields to more efficient  inference  than estimating the parameters in the models corresponding to the different groups separately. Moreover, the simultaneous estimation procedure can never be worse.
 While this property holds independently of the design under consideration, we  subsequently  construct efficient designs for  the comparison of curves  corresponding to
 not necessarily independent groups and demonstrate its superiority by means of a simulation study.

 The remaining part of this paper is organized as follows. In Section \ref{sec2} we introduce the basics  and the design problem.
 Section \ref{sec3} is devoted to a continuous time model, which could be interpreted as a limiting experiment of the   discrete model if the sample size converges to infinity. In this model we derive an explicit representation of the BLUE if estimation is performed simultaneously is both groups. In Section \ref{sec4} we develop a discrete approximation of the  continuous BLUE
by determining the optimal weights for the linear estimator. Finally, the optimal design points are determined such that the maximum width of the confidence band for the difference of the two regression functions is minimal. Section \ref{sec5} is devoted to a  small numerical comparison of the performance of the optimal designs with
uniform designs. In particular, it is demonstrated that
optimal designs yield substantially narrower confidence bands.   In many cases the maximal width  of a confidence band from the uniform design is by a factor between $ 2$ and $10$
larger than the width of a confidence band from the optimal
design.

\section{Simultaneous estimation of two regression models}
\label{sec2}
\def\theequation{2.\arabic{equation}}
\setcounter{equation}{0}
Throughout this paper we consider the
situation of two groups of observations
 $Y_{1,1} , \ldots , Y_{1,n}$
 and $Y_{2,1} , \ldots , Y_{2,n}$
at time points $t_1, \ldots , t_n$ ($i=1,2$)
where there may exist dependencies within and between the groups.
We assume the relation between the response and the covariate $t$ in  each group is described by a  linear regression models  given by
\begin{equation}\label{indivmod}
Y_{ij}= Y_i(t_j) = f_i^\top (t_j) \theta^{(i)} + {\eta}_i(t_j) \, , \, j=1, \ldots, n \, , i=1, 2 \, .
\end{equation}
Thus in each group  $n$ observations are taken at the same time points $t_1, \ldots, t_n$
which can be chosen in a compact interval, say $[a, b]$, and   observations at different time points and in different groups might be dependent.
  The vectors of the unknown parameters $\theta^{(1)}$ and $\theta^{(2)}$ are assumed to be $p_1$- and $p_2$-dimensional, respectively, and the corresponding vectors of regression functions  $f_i(t) = (f_{i,1}(t), \ldots, f_{i,p_i}(t))^\top$, $i=1, 2$, have  continuously differentiable
  linearly independent components.

To address the situation of correlation between the groups,
we start with a very simple covariance structure  for each group, but we emphasize that all results presented in this paper are correct for
more general covariance structures corresponding to Markov processes, see Remark \ref{rem1} for more details.
To be precise,  let
$\{{\varepsilon}_1(t)| ~t\in [a, b]\}$ and
$\{{\varepsilon}_2(t)| ~t\in [a, b]\}$
denote  two independent Brownian motions, such that
\begin{equation}\label{brown_props}
\mathbb{E}[\varepsilon_i(t_j)]= 0,~~K_i(t_j, t_k) = \mathbb{E}[\varepsilon_i(t_j)\varepsilon_i(t_k)] =
\min (t_j,t_k)
\end{equation}
denotes the mean value and the covariance of the individual process $\varepsilon_i$ at the points $t_{j}$ and $t_k$, respectively.  Let  $\sigma_1, \sigma_2 >0 $, $\varrho \in (-1, 1)$, denote by
$\mathbf\Sigma^{1/2}$   the square root of the covariance matrix
\be\label{Sigma}
\mathbf{\Sigma} = \begin{pmatrix} \sigma^2_1 & \sigma_1\sigma_2 \varrho \\ \sigma_1\sigma_2 \varrho & \sigma^2_2 \end{pmatrix} \,,
\ee
and  define for $t\in [a,b]$ the two-dimensional process
$ \{ \boldsymbol{\eta}(t)|~  t \in [a,b] \}  $ by
\begin{equation}\label{def_eta}
 \boldsymbol{\eta}(t)=
 \begin{pmatrix} \eta_1(t) \\  \eta_2 (t)
  \end{pmatrix}
 =\mathbf\Sigma^{1/2}
  \boldsymbol{\varepsilon}(t) ,
\end{equation}
where $  \boldsymbol{\varepsilon}(t) =(\varepsilon_1(t) , \varepsilon_2(t))^\top$.
Note that $\varrho \in (-1, 1)$ denotes the correlation between the observations  $Y_1(t_j)$ and $Y_2(t_j)$ ($j=1, \ldots , n$), and that in general the correlation between $Y_1(t_j)$ and $Y_2(t_k)$ is given by
$$
\mbox{Corr} (Y_1(t_j), Y_2(t_k)) = \varrho \left\{ \sqrt{\frac{t_j}{t_k}} \wedge \sqrt{\frac{t_k}{t_j}}\right\}
$$
if $t_j,t_k \in (0,1)$.

Considering  the two groups
individually   results in proper (for example weighted least squares) estimators
of the parameters $\theta^{(1)}$ and $\theta^{(2)}$. However, this procedure ignores the  correlation between the two groups and  estimating the parameters  $\theta^{(1)}$ and $\theta^{(2)}$
simultaneously from the
data of both groups might result in more precise estimates.
In order to define  estimators for  the
parameters  $\theta^{(1)}$ and $ \theta^{(2)}$
using the information from  both groups
we now consider a more general two-dimensional regression model,
which on the one hand contains the
 situation described in the previous paragraph as special case, but on the other hand
also allows us to consider  the  case, where some of the components in $\theta_1$ and $ \theta_2$ coincide, see 
Example \ref{ex2} and  Section
\ref{sec32} for details.
To be precise we define the  regression model
\be\label{moddisc}
\mathbf{Y}(t_j) = \begin{pmatrix} Y_1(t_j) \\ Y_2(t_j) \end{pmatrix} = \mathbf{F}^\top (t_j)\theta +
 \boldsymbol{\eta}(t_j)   =
 \mathbf{F}^\top (t_j)\theta + \mathbf{\Sigma}^{1/2}
 \boldsymbol{\varepsilon}(t_j), \quad \quad
    j=1, \ldots, n,
\ee
where   two-dimensional
observations
$$
\mathbf{Y}(t_1) = (Y_1(t_1), Y_2(t_1))^\top , \ldots ,  \mathbf{Y}(t_n) = (Y_1(t_n), Y_2(t_n))^\top
$$
   are  available
at   time points $t_1, \ldots, t_n \in [a,b]$.
In model \eqref{moddisc} the vector $\theta = (\vartheta_1, \ldots, \vartheta_p)^\top $
is a  $p$-dimensional parameter  and
\be \label{reg_mat}
\mathbf{F}^\top (t) = \begin{pmatrix} F^\top _1(t) \\ F^\top _2(t) \end{pmatrix} = \begin{pmatrix} F_{1,1}(t) & \ldots & F_{1,p}(t) \\ F_{2,1}(t) & \ldots & F_{2,p}(t) \end{pmatrix}
\ee
denotes a $(p\times 2)$ matrix
containing  continuously differentiable regression functions, where the two-dimensional functions
$(F_{1,1}(t) , F_{2,1}(t) )^\top ,  \ldots , (F_{1,p}(t) , F_{2,p}(t) )^\top$ are assumed to be linearly independent.

\begin{example} \label{ex1}
{\rm
The individual models defined in \eqref{indivmod} are contained in this two-dimensional model. More precisely,   defining the $p=(p_1+ p_2)$-dimensional vector of parameters $\theta$ by $\theta = ((\theta^{(1)})^\top , (\theta^{(2)})^\top )^\top $ and the regression function $\mathbf{F}^\top(t)$ in \eqref{reg_mat} by the rows
\begin{equation*}
F^\top _1(t)= ({f}^\top _1(t), 0^\top _{p_2}) , ~~ F^\top _2(t)= (0^\top _{p_1}, {f}^\top _2(t)),
\end{equation*}
it follows that   model \eqref{moddisc} coincides with model \eqref{indivmod}. Moreover, this composite model takes  the correlation between the groups into account. In this case the models describing the relation between the variable $t$ and the responses $Y_1(t)$ and $Y_2(t)$ do not share any parameters.
}
\end{example}

\begin{example} \label{ex2}
{\rm {
In this example we consider the situation where some of the parameters of the individual models in \eqref{indivmod} coincide. This situation occurs, for example,
if $Y_1(t)$  and $Y_2(t)$ represent clinical parameters 
(depending on time)
before and after treatment, where it can be assumed that the effect at time $a$ coincides before and after the treatment. 
In this case a reasonable model for average effect in 
the two groups is given by
$$
\mathbb{E} [Y_\ell (t) ] = \theta^{(0)} + (\tilde{\theta}^{(\ell )})^\top  \tilde {f}_\ell (t) 
~,~~\ell =1,2~.
$$
More generally, we consider the situation where the vectors  of the parameters are given by
$$\theta^{(1)}= (\theta^{(0)^\top }, \tilde{\theta}^{(1)^\top})^\top \quad , \quad \theta^{(2)}= (\theta^{(0)^\top}, \tilde{\theta}^{(2)^\top})^\top \, ,$$
where $\theta^{(0)} \in \mathbb{R}^{p_0}$ denotes the  vector of common parameters in both models
and  vectors $\tilde{\theta}^{(1)} \in \mathbb{R}^{p_1 - p_0}$  and $\tilde{\theta}^{(2)} \in \mathbb{R}^{p_2 - p_0}$ contain the different  parameters in  the two individual models. The corresponding regression functions are given by
\begin{equation}\label{holA}
    f_1^\top (t) = (f^\top_{0}(t), {\tilde{f}_1}^\top(t)) \quad, \quad f_2^\top (t) = (f^\top_{0}(t), \tilde{f}_2^\top(t)) \, ,
\end{equation}
where the vector $f^\top_0(t)$ contains the regression functions
corresponding to the common parameters in the two models, and  ${\tilde{f}_1}^\top(t)$ and  ${\tilde{f}_2}^\top(t)$   denote the
vectors of regression functions
corresponding to the different parameters $\tilde{\theta}^{(1)}$  and $\tilde{\theta}^{(2)}$,  respectively. \\
Defining the $p=(p_1+p_2-p_0)$-dimensional vector of parameters $\theta$ by $\theta= (\theta^{(0)}, \tilde\theta^{(1)}, \tilde\theta^{(2)})$ and the regression function $\mathbf{F}^\top(t)$ in \eqref{reg_mat} by the rows
\begin{equation*} 
F^\top _1(t)= (f^\top_0(t), \tilde{f}^\top _1(t), 0^\top _{p_2-p_0}) , ~~ F^\top _2(t)= (f^\top_0(t), 0^\top _{p_1-p_0}, \tilde{f}^\top _2(t)),
\end{equation*}
it follows that the model \eqref{moddisc} contains the individual models in \eqref{indivmod},
where the regression functions are given by \eqref{holA} and
the parameters $\theta^{(1)}$ and $\theta^{(2)}$ share the parameter $\theta^{(0)}$. Moreover,
this composite model takes
the potential correlation between the groups into account.
}
}
\end{example}



\section{Continuous time models}\label{sec3}
\def\theequation{3.\arabic{equation}}
\setcounter{equation}{0}

It was demonstrated by \cite{detkonzhi2017}  that
efficient designs for dependent data in regression problems can be derived by first  considering the estimation problem
in a continuous time model. In this model there is no optimal design problem as the data can be observed over the full interval $[a,b]$. However,
efficient designs can be determined in two steps.
First, one derives the best linear unbiased  estimator (BLUE)  in the continuous time model and, secondly,  one determines design points (and an estimator) such that the resulting estimator from the discrete data
provides a good   approximation of the optimal solution in the continuous time model.
In this  paper we will use this strategy
to develop optimal designs for the comparison of regression curves from  two (possible) dependent groups. In the present
section we discuss a continuous time model corresponding to the discrete model \eqref{moddisc}, while the second step, the determination of an   ``optimal'' approximation will be postponed to   the following  Section \ref{sec4} .

\subsection{Best linear unbiased estimation}
\label{sec30}

To be precise, we consider the  continuous time version of the linear regression model in \eqref{moddisc}, that is,
\be \label{mod1cont}
\mathbf{Y} (t)= \begin{pmatrix} Y_1(t) \\ Y_2(t) \end{pmatrix} = \mathbf{F}^\top (t)\theta + \mathbf{\Sigma}^{1/2}
 \boldsymbol{\varepsilon}(t)   \, ,  \quad t\in[a,b] \, ,
\ee
where we assume $0<a<b$ and the full trajectory of the process $\{ \boldsymbol Y(t) \mid t\in [a,b]\}$ is observed,  $\{\boldsymbol\varepsilon(t)=(\varepsilon_1(t), \varepsilon_2(t))^\top  \mid t\in [a,b]\}$ is a vector of independent Brownian motions as defined in \eqref{brown_props} and the matrix $\mathbf\Sigma^{1/2}$ is the square root of the covariance matrix $\mathbf{\Sigma}$ defined in \eqref{Sigma}.
Note that we restrict ourselves to an interval on the positive line, because in this case the notation is slightly simpler.   But we emphasize that the theory  developed in this section can also be applied for $a=0$, see Remark \ref{rem2} for more details.
We further assume that the $(p \times p)$-matrix 
\begin{equation}\label{cmatrix}
\mathbf{M} =\int_a^b \mathbf{\dot{F}}(t) \mathbf{\Sigma}^{-1} \mathbf{\dot{F}}^\top (t) \,dt + \frac{1}{a} \mathbf{{F}}(a)\mathbf{\Sigma}^{-1} \mathbf{{F}}^\top (a)
\end{equation}
is non-singular.\\

\begin{theorem}
\label{thm1}
Consider the continuous time linear regression model \eqref{mod1cont} on the interval $[a,b]$, $a >0$,
with a continuously differentiable matrix of regression
functions $\mathbf{F}$,  a vector
$\{\boldsymbol\varepsilon(t)=(\varepsilon_1(t), \varepsilon_2(t))^\top \mid t\in [a,b]\}$ of independent Brownian motions
and a covariance matrix $\mathbf{\Sigma}$ defined by \eqref{Sigma}.   The best linear unbiased estimator of the parameter $\theta$ is given by
\begin{equation}\label{asblue}
\hat{\theta}_{\rm BLUE}   = \mathbf{M}^{-1} \Big( \int_a^b \dot{\mathbf{F}}(t) \mathbf{\Sigma}^{-1}\,d\mathbf{Y} (t) + \frac{1}{a}  \mathbf{F}(a)\mathbf{\Sigma}^{-1} \mathbf{Y}(a) \Big) .
\end{equation}
Moreover, the minimum variance is given by
\begin{equation}\label{cmatrix_inv}
\mbox{\rm Cov}(\hat{\theta}_{\rm BLUE} ) =
\mathbf{M}^{-1} =\left(\int_a^b \mathbf{\dot{F}}(t) \mathbf{\Sigma}^{-1} \mathbf{\dot{F}}^\top (t) \,dt + \frac{1}{a} \mathbf{{F}}(a)\mathbf{\Sigma}^{-1} \mathbf{{F}}^\top (a)\right)^{-1} \, .
\end{equation}
\end{theorem}
\textbf{Proof.}
Multiplying
$ \boldsymbol{{Y}}$ by the matrix $\mathbf{\Sigma}^{-1/2}$
yields a transformed
regression model
\begin{equation}\label{contmod_uncorrelated}
\boldsymbol{\tilde{Y}}(t) = \begin{pmatrix}\tilde{Y}_1(t) \\ \tilde{Y}_2(t) \end{pmatrix}  = \mathbf{\Sigma}^{-1/2} \begin{pmatrix}{Y}_1(t) \\ {Y}_2(t) \end{pmatrix} = \mathbf{\Sigma}^{-1/2}\mathbf{F}^\top (t) \theta + \ \boldsymbol{\varepsilon}(t)  \, ,
\end{equation}
where
$\mathbf{\Sigma}^{-1/2}$ is the inverse of $\mathbf{\Sigma}^{1/2}$, the square root of the covariance matrix $\mathbf{\Sigma}$
defined in \eqref{Sigma}.
Note that the components of the vector  $\boldsymbol{\tilde{Y}}$ are independent, and consequently, the joint likelihood function can be obtained as  the product of the individual components. Next we rewrite the components of the continuous time model  \eqref{contmod_uncorrelated}
in terms  of two stochastic differential equations, that is
\begin{eqnarray}
d\tilde{Y}_i(t)& =& \mathbf{1}_{[a,b]}(t) \mathbf{\Sigma}^{-1/2}_i {\mathbf{\dot{F}}}^\top (t) \theta dt +  d\varepsilon_i(t) \, , t\in [0, b]~, \label{stochdgl} \\
\tilde{Y}_i(a) &= & \mathbf{\Sigma}^{-1/2}_i \mathbf{F}^\top (a) \theta + \varepsilon_i(a)~, \label{anfangswert}
\end{eqnarray}
 where $ \mathbf{1}_{A}$ is  the indicator function of the set $A$ and
$\mathbf{\Sigma}^{-1/2}_i$ denotes the $i$-th row of  the matrix $\mathbf{\Sigma}^{-1/2}$ ($i=1, 2$).
Since $\{\varepsilon_i(t)| ~t\in [a, b]\}$ is a Brownian motion its increments are independent. Consequently, the processes $\{\tilde{Y}_i(t)| ~t\in [0,b]\}$ and the random variable $\tilde Y_i(a)$ are independent. To obtain the joint density of the processes defined by \eqref{stochdgl} and \eqref{anfangswert} it is therefore sufficient  to derive the individual densities.
\\
Let $\mathbb{P}_\theta^{(i)}$  and $\mathbb{P}_0^{(i)}$  denote the measures on $C([0,b])$ associated with the process
 $\tilde Y_i = \{ Y_i(t) | \ t \in [0,b] \}$ and $\{ \varepsilon_{i} (t) | \ t \in [0,b] \}$, respectively.
It follows from Theorem 1 in Appendix II of
\cite{MR620321}
 that $\mathbb{P}_\theta^{(i)} $ is absolute continuous with respect to $\mathbb{P}_0^{(i)}$
 with
 Radon-Nikodym-density
\begin{equation*}
\frac{d\mathbb{P}_\theta^{(i)}} 
{d \mathbb{P}_0^{(i)}} ( \tilde{Y}_i) = \exp\left\{\int_a^b\Sigma^{-1/2}_i \mathbf{\dot{F}}^\top(t)\theta d\tilde{Y}_i(t) - \frac{1}{2} \int_a^b (\Sigma^{-1/2}_i \mathbf{\dot{F}}^\top(t)\theta)^2 dt\right \}\, .
\end{equation*}
Similarly, if $\mathbb{Q}_\theta $ denotes the distribution of the random variable $\tilde Y_i(a) \sim {\cal N } (\Sigma^{-1/2}_i \mathbf{F}^\top (a)\theta, a)
$ in
\eqref{anfangswert}, then
the Radon-Nikodym-density of
$\mathbb{Q}^{(i)}_\theta $ with respct to $\mathbb{Q}^{(i)}_0 $
is given by
\begin{equation*}
    \frac{d\mathbb{Q}^{(i)}_\theta}{d \mathbb{Q}^{(i)}_0}
    (\tilde{Y}_i(a))= \exp\left\{\frac{\tilde{Y}_i(a)\Sigma^{-1/2}_i \mathbf{F}^\top (a)\theta}{a} - \frac{1}{2} \frac{(\mathbf{\Sigma}^{-1/2}_i \mathbf{F}^\top (a) \theta)^2}{a} \right\} \, .
\end{equation*}
Consequently, because of independence,
the joint density of $(\mathbb{P}^{(i)}_\theta,\mathbb{Q}^{(i)}_\theta)$
with respect to $(\mathbb{P}^{(i)}_0,\mathbb{Q}^{(i)}_0)$
is obtained as
\begin{equation*}
\begin{split}
\frac{d\mathbb{P}^{(i)}_\theta}{d\mathbb{P}^{(i)}_0} (\tilde{Y}_i)\times \frac{d\mathbb{Q}^{(i)}_\theta}{d\mathbb{Q}^{(i)}_0}
(\tilde{Y}_i(a)) =&  \exp\left\{\left(\int_a^b\mathbf{\Sigma}^{-1/2}_i \mathbf{\dot{F}}^\top (t)\theta d\tilde{Y}_i(t)+ \frac{\tilde{Y}_i(a)\mathbf{\Sigma}^{-1/2}_i \mathbf{F}^\top (a)\theta}{a} \right) \right.  \\
   & \left. - \frac{1}{2}\left( \int_a^b (\mathbf{\Sigma}^{-1/2}_i \mathbf{\dot{F}}^\top (t)\theta)^2 dt\ + \frac{(\mathbf{\Sigma}^{-1/2}_i \mathbf{F}(a) \theta)^2}{a}\right) \right\} \, .
\end{split}
\end{equation*}
As the processes $\tilde{Y}_1$ and $\tilde{Y}_2$
are independent by construction the maximum likelihood estimator in the model \eqref{mod1cont}
can be determined by  solving  the  equation
\begin{equation*}
\begin{split}
&  \frac{\partial}{\partial \theta} \log \Big \{ \prod_{i=1}^{2} \frac{d\mathbb{P}^{(i)}_\theta}{d\mathbb{P}^{(i)}_0} (\tilde{Y}_i)\times \frac{d\mathbb{Q}^{(i)}_\theta}{d\mathbb{Q}^{(i)}_0} (\tilde{Y}_i(a))\Big  \}  =
\sum_{i=1}^{2} \Big  \{ \int_a^b\mathbf{\dot{F}}(t)\mathbf{\Sigma}^{-1/2}_i  d\tilde{Y}_i(t)+ \frac{\mathbf{F}(a)\mathbf{\Sigma}^{-1/2}_i\tilde{Y}_i(a)}{a} \\
 & ~~~~~~  ~~~~~~  ~~~~~~   -   \Big( \int_a^b \mathbf{\dot{F}}(t)\Sigma^{-1/2}_i \mathbf{\Sigma}^{-1/2}_i  \mathbf{\dot{F}}^\top (t) \,dt + \mathbf{\dot{F}}(a)\mathbf{\Sigma}^{-1/2}_i\mathbf{\Sigma}^{-1/2}_i \mathbf{\dot{F}}^\top (a) \Big) \theta
\Big  \}= 0
\end{split}
\end{equation*}
with respect to $\theta$.
The solution coincides with  the linear estimate defined in \eqref{asblue}, and a straightforward calculation, using Ito's formula
and the fact that the random variables $\int^b_a \boldsymbol{\dot{F}} (t) d \boldsymbol{\varepsilon} _t $
and $ \boldsymbol{\varepsilon}_a$ are independent, gives
\begin{eqnarray*}
 \mbox{\rm Cov}(\hat{\theta}_{\rm BLUE} ) &=& \mathbf{M}^{-1} \mathbb{E}_{\theta} \Bigl [
  \Big( \int_a^b \dot{\mathbf{F}}(t) \mathbf{\Sigma}^{-1}\,d\mathbf{Y}(t) + \frac{1}{a}  \mathbf{F}(a)\mathbf{\Sigma}^{-1} \mathbf{Y}(a) \Big) \\
 && \times
 \Big( \int_a^b \dot{\mathbf{F}}(t) \mathbf{\Sigma}^{-1}\,\mathbf{Y}(t) + \frac{1}{a}  \mathbf{F}(a)\mathbf{\Sigma}^{-1} \mathbf{Y}(a)  \Big)^\top  \Bigr] \mathbf{M}^{-1} \\
&=& \mathbf{M}^{-1} \Big ( \int_a^b \mathbf{\dot{F}}(t) \mathbf{\Sigma}^{-1} \mathbf{\dot{F}}^\top (t) \,dt + \frac{1}{a} \mathbf{{F}}(a)\mathbf{\Sigma}^{-1} \mathbf{{F}}^\top (a)\Big) \mathbf{M}^{-1} = \mathbf{M}^{-1},
  \end{eqnarray*}
where the matrix $\mathbf{M}$ is defined in \eqref{cmatrix}.
Since the covariance matrix $\mathbf{M}^{-1}$ is the inverse of the information matrix in   the continuous time regression model in \eqref{mod1cont} \citep[see][p. 81]{MR620321}, the linear estimator \eqref{asblue} is the BLUE, which completes the proof of  Theorem \ref{thm1}.
\hfill $\Box$

\medskip 

\begin{remark} \label{rem2}
 {\rm ~
The proof of Theorem \ref{thm1} can easily be modified to obtain the BLUE for the continuous time model on the interval   $a=[0, b]$. More precisely, for $a=0$   equation   \eqref{anfangswert} becomes a deterministic equation equivalent to
\begin{equation}
\mathbf{Y}(0) = \mathbf{F}^\top (0) \theta \, \label{anfangswert_a0} \, , 
\end{equation}
and we  have to distinguish three cases. 
\begin{enumerate}
    \item[(1)] If the regression function $\mathbf{F}$ satisfies $\mathbf{F}(0) = \mathbf{0}_{p\times 2}$ (that is $\mbox{rank}(\mathbf{F}(0)) = 0$)), the deterministic equation   \eqref{anfangswert_a0} does not  contain any further information about the parameter $\theta$ and the maximum likelihood estimator in model \eqref{mod1cont} is given by
\begin{equation*}\label{contblue_a0}
\hat{\theta}_{\rm BLUE}   = \mathbf{M}_0^{-1} \Big( \int_0^b \dot{\mathbf{F}}(t) \mathbf{\Sigma}^{-1}\,d\mathbf{Y} (t) \Big) \, , 
\end{equation*}
where the minimum variance is given by
\begin{equation*}\label{m0matrix_inv}
\mbox{\rm Cov}(\hat{\theta}_{\rm BLUE} ) =
\mathbf{M}_0^{-1} =\left(\int_0^b \mathbf{\dot{F}}(t) \mathbf{\Sigma}^{-1} \mathbf{\dot{F}}^\top (t) \,dt\right)^{-1} \, .
\end{equation*}
\item[(2)] If the rank of the matrix $\mathbf{F}(0)$ satisfies $\mbox{rank}(\mathbf{F}(0))=1$, the deterministic equation \eqref{anfangswert_a0} contains one informative equation about $\theta$. In that case, we assume without loss of generality that $F_{1,1}(0) \neq 0 $ and it follows by \eqref{anfangswert_a0} that $\theta_1$ can be reformulated by $\theta_2, \ldots, \theta_p$ through 
\begin{equation}\label{theta_red_case2}
    \theta_1 = \frac{Y_1(0) - \sum_{i=j}^p\theta_j F_{1,j}(0)}{F_{1,1}(0)} \, .
\end{equation}
Using \eqref{theta_red_case2} in combination with model \eqref{mod1cont}, we obtain a reduced model by
\begin{equation}\label{modreduced_case2}
\mathbf{Z}(t) = \mathbf{Y}(t) - \frac{Y_1(0)}{F_{1,1}(0)}\begin{pmatrix}F_{1, 1}(0) \\ F_{2, 1}(0)\end{pmatrix} =  \tilde{\mathbf{F}}(t)\tilde{\theta}  + \mathbf{\Sigma}^{1/2}
 \boldsymbol{\varepsilon}(t) \, ,
\end{equation}
where the matrix valued function $\tilde{\mathbf{F}}(t)$ is defined by
\begin{equation}\label{reduced_fun_case2}
    \tilde{\mathbf{F}}^T(t) = \bigl( F_{i,j}(t) - \frac{F_{i,1}(0)}{F_{1, 1}(0)} F_{1, j}(0)\bigr)_{i=1, 2, j=2, \ldots p}
\end{equation}
and the reduced $(p-1)$-dimensional parameter $\tilde\theta$ is given by $\tilde\theta = (\theta_2, \ldots, \theta_p) $ . 
It follows by $\mbox{rank}(\mathbf{F}(0)) = 1$, that the matrix valued function $\tilde{\mathbf{F}}(t)$ defined in \eqref{reduced_fun_case2} satisfies $\mathbf{\tilde{F}^T(0)} = 0_{2\times p}$. Consequently, the modified model given by \eqref{modreduced_case2} satisfies the condition of the case given in (1)  and the best linear unbiased estimator for the reduced parameter $\tilde\theta$ is obtained by 
 \begin{equation}\label{contblue_case2}
\hat{\tilde\theta}_{\rm BLUE}   = \mathbf{M}_0^{-1} \Big( \int_0^b \dot{\tilde{\mathbf{F}}}(t) \mathbf{\Sigma}^{-1}\,d\mathbf{Z} (t) \Big) \, , 
\end{equation}
where the process $\{\mathbf{Z(t)}; t\in [0, b]\}$ is defined by \eqref{modreduced_case2}, the matrix $\mathbf{\tilde{F}}(t)$ is given by \eqref{reduced_fun_case2} and the minimum variance is given by
\begin{equation*}\label{m0matrix_inv_case2}
\mbox{\rm Cov}(\hat{\tilde\theta}_{\rm BLUE}  ) =
\mathbf{M}_0^{-1} =\left(\int_0^b \dot{\tilde{\mathbf{F}}}(t) \mathbf{\Sigma}^{-1} \dot{\tilde{\mathbf{F}}}^\top (t) \,dt\right)^{-1} \, .
\end{equation*}
The best linear unbiased estimator for the remaining parameter $\theta_1$ is then obtained by 
$$\hat\theta_1 = \frac{Y_1(0) - \sum_{i=j}^p\hat{\tilde\theta}_{{\rm BLUE},j} F_{1,j}(0)}{F_{1,1}(0)} \, .$$
\item[(3)] If the rank of the matrix $\mathbf{F}(0)$ satisfies $\mbox{rank}(\mathbf{F}(0)) = 2$, equation \eqref{anfangswert_a0} contains two informative equations about $\theta$. \\
Let 
\begin{equation}\label{Amat}
    \mathbf{A}(t) = \begin{pmatrix} F_{1,1}(t) & F_{1,2}(t) \\ F_{2,1}(t) & F_{2, 2}(t) \end{pmatrix}
\end{equation}
be the submatrix of $\mathbf{F}$ which contains the first two columns of $\mathbf{F}^T(t)$. Without loss of generality, we assume that $\mathbf{A}(0)$ is non-singular (as $\mbox{rank}(\mathbf{F}(0)) = 2$). \\
Then it follows by \eqref{anfangswert_a0} that 
\begin{equation}\label{reduced_theta_case3}
    \begin{pmatrix}
    \theta_1 \\ \theta_2 
    \end{pmatrix}  = \mathbf{A}^{-1}(0)\Bigl(\mathbf{Y}(0) - \bigl(\sum_{j=3}^{p} F_{i,j}(0)\theta_j\bigr)_{i=1, 2} \Bigr) \, .
\end{equation}
Using \eqref{reduced_theta_case3} in combination with \eqref{mod1cont} we obtain a reduced model given by
\begin{equation}\label{reducedmod_case3}
\mathbf{Z}(t) = \mathbf{Y}(t) - \mathbf{A}(t)\mathbf{A}^{-1}(0) \mathbf{Y}(0) = \tilde{\mathbf{F}}(t)\tilde{\theta}  + \mathbf{\Sigma}^{1/2}
 \boldsymbol{\varepsilon}(t)
\end{equation}
where the matrix valued function $\mathbf{A}(t)$ is given by \eqref{Amat}, the matrix valued function $\tilde{\mathbf{F}}^T(t)$ is of the form
\begin{equation}\label{reduced_fun_case3}
    \tilde{\mathbf{F}}^T(t) = \mathbf{A}(t) \mathbf{A}^{-1}(0)\bigl(F_{i,j}(0)\bigr)_{i=1, 2; j=3, \ldots p } + \bigl(F_{i,j}(t)\bigr)_{i=1, 2; j=3, \ldots p }
\end{equation}
and the reduced $(p-2)$-dimensional parameter $\tilde\theta$ is given by $\tilde\theta = (\theta_3, \ldots, \theta_p) $ .  \\
The matrix valued function $\tilde{\mathbf{F}}(t)$ defined in \eqref{reduced_fun_case3} satisfies $\mathbf{\tilde{F}^T}(0) = \mathbf{0}_{2\times p}$. 
 Consequently, the modified model given by \eqref{reducedmod_case3} satisfies the condition of the case given in (1) and the best linear unbiased estimator $\hat{\tilde\theta}_{\rm BLUE}$ for the reduced $(p-2)$-dimensional parameter $\tilde\theta$ is obtained by \eqref{contblue_case2} using the process $\{\mathbf{Z(t)}; t\in [0, b]\}$ defined by \eqref{reducedmod_case3} and the matrix valued function $\mathbf{\tilde{F}}(t)$ given by \eqref{reduced_fun_case3}. 
 The best linear unbiased estimator for the remaining parameter $(\theta_1, \theta_2)^T$ is then obtained by 
$$\begin{pmatrix}
    \hat\theta_1 \\ \hat\theta_2 
    \end{pmatrix}  = \mathbf{A}^{-1}(0)\Bigl(\mathbf{Y}(0) - \bigl(\sum_{j=3}^{p} F_{i,j}(0)\hat{\tilde\theta}_{\rm BLUE,j}\bigr)_{i=1, 2} \Bigr) \, .$$

\end{enumerate}

}
\end{remark}

\subsection{Model with no common parameters}
\label{sec31}
Recall the definition of   model \eqref{indivmod}
in Section \ref{sec1}.
It was demonstrated in Example \ref{ex1} that this case
is a special case of   model \eqref{moddisc},
where the matrix  $\mathbf{F}^\top$ is given by
\be \label{modcont_sep1}
\mathbf{F}^\top (t)
=
\begin{pmatrix} {f}^\top _1(t) & 0^\top _{p_2}\\ 0^\top _{p_1} &  {f}^\top _2(t)\end{pmatrix}
\ee
and $\theta  =({\theta^{(1)}}^\top ,{\theta^{(2)}}^\top )^\top$.
Considering both components in the vector $\mathbf{Y}$ separately, we obtain
continuous versions of the two models
introduced in  \eqref{indivmod},that is,
\begin{equation}\label{indiv_contmod}
Y_i(t) = {f}^\top _i(t) \theta^{(i)} + \eta_i(t), \, i=1, 2 \, ,
\end{equation}
where the error processes $\{\eta(t) \mid t \in [a, b]\}$ is defined by \eqref{def_eta}.
An application of  Theorem \ref{thm1} yields the following
BLUE.
\medskip

\begin{corollary}\label{est_noncommon}
Consider the continuous time linear regression model \eqref{moddisc} on the interval $[a, b]$, with
continuously differentiable  matrix \eqref{modcont_sep1},  a vector $\{\boldsymbol\varepsilon(t)=(\varepsilon_1(t), \varepsilon_2(t))^\top \mid t\in [a,b]\}$
 of independent Brownian motions and a  matrix $\mathbf{\Sigma}$ defined by \eqref{Sigma}.   The best linear unbiased estimator for the parameter $\theta$ is given by
\begin{equation} \label{bluesep}
\begin{split}
    \hat\theta_{\rm BLUE} = \begin{pmatrix} \hat\theta^{(1)}_{\rm BLUE}\\\hat\theta^{(2)}_{\rm BLUE} \end{pmatrix}= \frac{1}{\sigma^2_1\sigma^2_2(1-\varrho^2)} \mathbf{M} ^{-1} \left\{ \int_a^b \begin{pmatrix}
\sigma_2^2 \dot{f}_1(t) & - \sigma_1\sigma_2\varrho \dot{f}_1(t) \\
  - \sigma_1\sigma_2\varrho  \dot{f}_2(t) &  \sigma^2_1\dot{f}_2(t)
\end{pmatrix}
d\begin{pmatrix}Y_1(t) \\ Y_2(t) \end{pmatrix} \, \right. \\
\left. + \frac{1}{a} \begin{pmatrix}
\sigma_2^2 {f}_1(a) & - \sigma_1\sigma_2\varrho {f}_1(a) \\
  - \sigma_1\sigma_2\varrho  {f}_2(a) &  \sigma^2_1{f}_2(a)
  \end{pmatrix}\begin{pmatrix}Y_1(a) \\ Y_2(a) \end{pmatrix}
  \right\} \, .
  \end{split}
\end{equation}
The minimum variance is given by  $\mathbf{M}^{-1}$, where
\begin{equation*} 
\mathbf{M} = \frac{1}{\sigma^2_1\sigma^2_2(1-\varrho^2)} \begin{pmatrix} \sigma_2^2 \mathbf{M}_{11} & -\sigma_1\sigma_2\varrho \mathbf{M}_{12} \\ -\sigma_1\sigma_2\varrho\mathbf{M}_{21} & \sigma^2_1\mathbf{M}_{22}\end{pmatrix}
\end{equation*}
and
\begin{equation}\label{indiv_mat}
\mathbf{M}_{ij} = \int_a^b  \dot{f}_i(t) \dot{f}^\top _j(t)dt \, + \frac{1}{a} f_i(a) f_j^T(a) \, ~~~~i,j=1, 2 .
\end{equation}
\end{corollary}


It is of interest to compare the estimator \eqref{bluesep} with the estimator $\hat\theta_{\rm mar} = ((\hat\theta_{\rm mar}^{(1)})^\top, (\hat\theta_{\rm mar}^{(2)})^\top)^\top $, which is obtained by estimating the parameter in    both models \eqref{indiv_contmod} separately. It follows from Theorem 2.1 in \cite{detkonzhi2017} that the best linear unbiased estimators
in these models are given by
\begin{equation}\label{individual_est}
\hat{\theta}^{(\ell)}_{\rm mar} =\mathbf{M}^{-1}_{\ell\ell}\left( \int_{a}^b \dot{{f}}_\ell(t) dY_\ell(t) + \frac{1}{a}{f}_\ell(a) Y_\ell(a) \right) \, , \quad \quad
\ell =1,2 ,
\end{equation}
where the matrices
are defined by
$$
\mathbf{M}_{\ell\ell} = \int_a^b  \dot{f}_\ell(t) \dot{f}^\top _\ell(t)dt + \frac{1}{a}f_\ell(a) f_\ell^\top(a)
, ~~~\ell=1, 2.
$$
Moreover, the covariance matrices of the  estimators $\hat\theta^{(1)}_{\rm mar}$ and $\hat\theta^{(2)}_{\rm mar}$ are the inverses of the Fisher information matrices in the individual models, that
is
\begin{equation} \label{hol1}
    \mbox{Cov}(\hat\theta^{(\ell)}_{\rm mar}) = \sigma_\ell^2 \mathbf{M}^{-1}_{\ell\ell}
~~~\ell=1,2.
\end{equation}
The following result compares the variance of the two estimators
\eqref{bluesep} and \eqref{individual_est}.

\begin{theorem}\label{loewner_individ_mod}
If the assumptions of Corollary \ref{est_noncommon} are satisfied, we have (with respect to the Loewner ordering)
$$
\mbox{\rm Cov}(  \hat\theta^{(\ell )}_{\rm BLUE} )
 \leq
 \mbox{\rm Cov}(\hat\theta^{(\ell)}_{\rm mar})
~, ~~ \ell =1,2 \, , 
$$
for all $\varrho \in (-1, 1)$, where the
$ \hat\theta^{(\ell )}_{\rm BLUE} $ and $ \hat\theta^{(\ell)}_{\rm mar}$ are the best linear unbiased estimators of the parameter $\theta^{(\ell)}$
obtained by simultaneous estimation (see \eqref{bluesep})
and separate estimation  in the two groups (see \eqref{individual_est}) ,
respectively.
\end{theorem}
\textbf{Proof.}
Without loss of generality we consider the case $\ell =1$, the proof for the index $\ell=2$ is obtained by the same arguments.
Let $\mathbf{K_1}^\top  = (\mathbf{I}_{p_1}, \mathbf{0}_{p_1\times p_2})$ be a $p_1\times (p_1+ p_2)$- matrix,
where $\mathbf{I}_{p_1}$ and $\mathbf{0}_{p_1\times p_2}$ denote the $p_1$-identity matrix and a $(p_1\times p_2)$-matrix filled with zeros. Then,
$$
\mbox{\rm Cov}(  \hat\theta^{(\ell )}_{\rm BLUE} )
= (\mathbf{C}_{\mathbf{K}_1}(\mathbf{M}))^{-1}, $$
where
\begin{equation}\label{schurkomplement}
\mathbf{C}_{\mathbf{K}_1}(\mathbf{M}) = (\mathbf{K}^\top _1\mathbf{M}^{-1}\mathbf{K}_1)^{-1}
= \frac{1}{\sigma^2_1(1-\varrho^2)} \left( \mathbf{M}_{11} -\varrho^2 \mathbf{M}_{12}\mathbf{M}^{-1}_{22}\mathbf{M}^\top _{12}\right)
\end{equation}
is the  Schur  complement of the block $\mathbf{M}_{22}$  of the information matrix $\mathbf{M}$ \citep[see p.~74 in][]{pukelsheim2006}.
Observing \eqref{hol1} we  now compare $\mathbf{C}_{\mathbf{K}_1}(\mathbf{M})$ and $\tfrac{1}{\sigma^2}\mathbf{M}_{11}$ and obtain
\begin{equation}\label{umformung_infos}
\begin{split}
\mathbf{C}_{\mathbf{K}_1}(\mathbf{M})-\frac{1}{\sigma^2_1}\mathbf{M}_{11} =& \frac{1}{\sigma^2_1(1-\varrho^2)} \left( \mathbf{M}_{11} - \varrho^2 \mathbf{M}_{12}\mathbf{M}^{-1}_{22} \mathbf{M}_{12}^\top \right) - \frac{1}{\sigma^2_1}\mathbf{M}_{11} \\
=& \frac{\varrho^2}{\sigma^2_1(1-\varrho^2)}\left( \mathbf{M}_{11} - \mathbf{M}_{12}\mathbf{M}^{-1}_{22} \mathbf{M}_{12}^\top \right) \,  \\
:=& \frac{\varrho^2}{\sigma^2_1(1-\varrho^2)} \mathbf{C}_{\mathbf{K}_1}({\mathbf{\tilde M}})\, ,
\end{split}
\end{equation}
where  $\mathbf{C}_{\mathbf{K}_1}({\mathbf{\tilde M}})$ is the Schur complement of the
block $\mathbf{M}_{22}$ of the matrix
\begin{equation*} 
\mathbf{\tilde M}  = \begin{pmatrix} \mathbf{M}_{11} & \mathbf{M}_{12} \\ \mathbf{M}_{21} & \mathbf{M}_{22}\end{pmatrix} \, .
\end{equation*}
Note that the matrix $\mathbf{\tilde M}$ is nonnegative definite. An application of  Lemma 3.12 of \cite{pukelsheim2006} shows
that the Schur complement $\mathbf{C}_{\mathbf{K}_1}(\mathbf{\tilde M})$ is also nonnegative definite, that is $\mathbf{C}_{\mathbf{K}_1}(\mathbf{\tilde M})\geq 0 $ with respect to the Loewner ordering. Observing  \eqref{umformung_infos}  we have
$$
\big (
\mbox{\rm Cov}(  \hat\theta^{(1 )}_{\rm BLUE} )\big )^{-1}
=
 \mathbf{C}_{\mathbf{K}_1}(\mathbf{M}) \geq   \frac{1}{\sigma^2_1}\mathbf{M}_{11}
=
\big (  \mbox{\rm Cov}(\hat\theta^{(1)}_{\rm mar}) \big )^{-1}
 $$
and
the statement of the theorem follows.
 \hfill $\Box$
\bigskip

\begin{remark} \label{vergleich}
{\rm
    If  $\varrho = 0$ we have  $\mathbf{C}_{K_1}(\mathbf{M})= \mathbf{M}_{11}$, and it follows
    from
    \eqref{schurkomplement}
  that separate  estimation in the individual groups does not yield less precise estimates, that is
  $
\mbox{\rm Cov}(\hat\theta_{\rm mar}^{(l)})
= \mbox{\rm Cov}(\hat\theta^{(1)}_{\rm BLUE})$
$(\ell =1,2)$.
  However,  in general we have $\mbox{\rm Cov} (\hat\theta_{\rm mar}^{(l)})
\geq \mbox{\rm Cov} (\hat\theta^{(1)}_{\rm BLUE})$.
Moreover,
the inequality is strict in most cases, which means that simultaneous estimation of the parameters $\theta^{(1)}$
and  $\theta^{(2)}$ yields more precise estimators.
A necessary condition for strict inequality
(i.e the matrix $\mbox{\rm Cov} (\hat\theta_{\rm mar}^{(l)})
 - \mbox{\rm Cov}(\hat\theta^{(1)}_{\rm BLUE})$
is positive definite) is the condition
$\varrho \not = 0$. The following  result
shows that this condition is not sufficient.
It considers the important case where the regression functions $f_1$ and $f_2$ in \eqref{modcont_sep1} are the same and shows that in this case the two estimators $\hat\theta_{\rm BLUE}$ and  $\hat\theta_{\rm mar}$ coincide.

}
\end{remark}

  \begin{corollary}
  If the assumptions of Corollary \ref{est_noncommon}  hold and additionally the regression functions in
  model   \eqref{moddisc}  satisfy $f_1 = f_2$,
  the best linear unbiased estimator for the parameter $\theta$  is given by
\begin{equation*}  
\begin{split}
    \hat\theta_{\rm BLUE} = \begin{pmatrix} \hat\theta^{(1)}_{\rm BLUE}\\\hat\theta^{(2)}_{\rm BLUE} \end{pmatrix}=
    \int_a^b \left(\mathbf{I}_2 \otimes \mathbf{M}_{11}^{-1}\dot{f}_1(t)\right) d\mathbf{Y} (t) + \frac{1}{a}\left(\mathbf{I}_2\otimes \mathbf{M}_{11}^{-1} f_1(a)\right) \mathbf{Y} (a)
  \end{split}\, ,
\end{equation*}
where $\mathbf{I}_{2}$ denotes the $2\times 2$-identity matrix and the matrix $\mathbf{M}_{11}$ is defined by \eqref{indiv_mat}.
Moreover,  the minimum variance is given by
$\mbox{\rm Cov}(\hat\theta_{\rm BLUE})= \mathbf{\Sigma} \otimes \mathbf{M}_{11}^{-1} $ and
$$
\mbox{\rm Cov}(\hat\theta_{\rm mar}^{(l)})
= \mbox{\rm Cov}(\hat\theta^{(1)}_{\rm BLUE}) ~~~(\ell =1,2)~.
$$
\end{corollary}

\subsection{Models with common parameters  }
\label{sec32}
 Recall the definition of  model \eqref{indivmod}
 in Section \ref{sec1}.
 It was demonstrated in Example \ref{ex2} that this case
 is a special case of  model \eqref{moddisc},
 where the matrix  of regression functions
 is given by
 \begin{equation}\label{commonF}
 \mathbf{F}^\top (t)
 =
 \begin{pmatrix} f^\top_0(t), \tilde{f}^\top _1(t), 0^\top _{p_2-p_0}\\f^\top_0(t), 0^\top _{p_1-p_0}, \tilde{f}^\top _2(t) \end{pmatrix}
 \end{equation}
 and the vector of  parameters is defined by
 $$\theta= (\theta^{(0)^\top}, \tilde\theta^{(1)^\top}, \tilde\theta^{(2)^\top})^\top \, . $$
 An application of  Theorem \ref{thm1} yields the
 BLUE in model \eqref{moddisc}
 with the matrix $ \mathbf{F}^\top$
 defined by \eqref{commonF}.

 \begin{corollary}\label{est_common}
 Consider the continuous time linear regression model
\eqref{moddisc}
 on the interval $[a, b]$,
 where the  the matrix  of regression functions $ \mathbf{F}^\top$
is  continuously differentiable.
The best linear unbiased estimator for the parameter $\theta$ is given by

 \begin{equation}\label{bluecom}
 \begin{split}
     \hat\theta_{\rm BLUE} = \begin{pmatrix} \hat\theta^{(0)}_{\rm BLUE} \\ \hat{\tilde\theta}^{(1)}_{\rm BLUE}\\\hat{\tilde\theta}^{(2)}_{\rm BLUE} \end{pmatrix}= \frac{1}{\sigma_1^2\sigma_2^2(1-\varrho^2)} \mathbf{M} ^{-1}
     \left\{\int_a^b \begin{pmatrix}
 (\sigma_2^2-\sigma_1\sigma_2\varrho)\dot{f}_0(t) & (\sigma_1^2-\sigma_1\sigma_2\varrho)\dot{f}_0(t)  \\
 \sigma_2^2\dot{\tilde{f}}_1(t)  & - \sigma_1\sigma_2 \varrho \dot{\tilde{f}}_1(t) \\
 - \sigma_1\sigma_2 \varrho \dot{\tilde{f}}_2(t) &  \sigma^2_1\dot{\tilde{f}}_2(t)
 \end{pmatrix}
 d\begin{pmatrix}Y_1(t) \\ Y_2(t) \end{pmatrix} \right. \\
 \left. +\frac{1}{a} \begin{pmatrix}
 (\sigma_2^2-\sigma_1\sigma_2\varrho){f}_0(a) & (\sigma_1^2-\sigma_1\sigma_2\varrho){f}_0(a)  \\
 \sigma_2^2{\tilde{f}}_1(a)  & - \sigma_1\sigma_2 \varrho {\tilde{f}}_1(a) \\
 - \sigma_1\sigma_2 \varrho {\tilde{f}}_2(a) &  \sigma^2_1{\tilde{f}}_2(a)
 \end{pmatrix}\,\begin{pmatrix}Y_1(a) \\ Y_2(a) \end{pmatrix} \right\} .
 \end{split}
 \end{equation}
The minimum variance is
 $$
 \mbox{\rm Cov}(  \hat\theta_{\rm BLUE} ) =
 \mathbf{M}^{-1}~,
 $$
where
 \begin{equation*}
 \mathbf{M} = \frac{1}{\sigma_1^2\sigma_2^2(1-\varrho^2)} \begin{pmatrix} (\sigma_1^2+\sigma_2^2 - \sigma_1\sigma_2\varrho)\mathbf{M}_{00} & (\sigma_2^2-\sigma_1\sigma_2\varrho)\mathbf{M}_{01} & (\sigma_1^2-\sigma_1\sigma_2\varrho)\mathbf{M}_{02}\\
 (\sigma_2^2-\sigma_1\sigma_2\varrho)\mathbf{M}_{10} & \sigma_2^2\mathbf{M}_{11} &  - \sigma_1\sigma_2 \varrho \mathbf{M}_{12} \\ (\sigma_1^2-\sigma_1\sigma_2\varrho)\mathbf{M}_{20} & - \sigma_1\sigma_2 \varrho\mathbf{M}_{21} & \sigma_1^2\mathbf{M}_{22}\end{pmatrix}
 \end{equation*}
and  individual blocks  in this matrix are given by
 \begin{equation}\label{indiv_mat}
 \mathbf{M}_{ij} = \int_a^b  \dot{g}_i(t) \dot{g}^\top _j(t)dt + \frac{1}{a}g_i(a) g_j^\top(a)\, ,
 \end{equation}
 for $i, j=0, 1, 2$, where $g_0(t) = f_0(t)$ and $g_i(t)= \tilde{f}_i(t)$ for $i=1, 2$ .
 \end{corollary}

 It is again of interest to compare the estimate \eqref{bluecom} with the estimate $\hat\theta_{\rm mar} = ((\hat\theta_{\rm mar}^{(1)})^\top, (\hat\theta_{\rm mar}^{(2)})^\top)^\top $, which is obtained by estimating the parameter $\theta^{(1)} = ((\theta^{(0)})^\top, (\tilde\theta^{(1)})^\top )^\top $ in both models \eqref{indiv_contmod} separately by using \eqref{individual_est}. The corresponding covariances of the estimators $\hat\theta^{(1)}_{{\rm mar}}$ and  $\hat\theta^{(2)}_{{\rm mar}}$ are given by \eqref{hol1}.
 The following result compares the variance of the two estimators \eqref{bluecom} and \eqref{individual_est}.
Its  proof is similar to the proof of Theorem \ref{loewner_individ_mod} and therefore omitted.
 \begin{theorem}\label{loewner_common}
If the assumptions of Corollary \ref{est_common} are satisfied, we have (with respect to the Loewner ordering)
$$
\mbox{\rm Cov}(  \hat\theta^{(\ell )}_{\rm BLUE} )
 \leq
 \mbox{\rm Cov}(\hat\theta^{(\ell)}_{\rm mar})
~, ~~ \ell =1,2 \, , 
$$
for all $\varrho \in (-1, 1)$, where the
$ \hat\theta^{(\ell )}_{\rm BLUE} $ and $ \hat\theta^{(\ell)}_{\rm mar}$ are the best linear unbiased estimators of the parameter $\theta^{(\ell)}$
obtained by simultaneous
and separate estimation,
respectively.
\end{theorem}

\begin{remark} \label{rem1}
{\rm ~
The results presented so far have been derived  for the case where the error process $\{\boldsymbol{\varepsilon}(t) = (\varepsilon_1(t), \varepsilon_2(t))^\top| ~t\in[a, b]\}$
in \eqref{moddisc}
consist of two independent Brownian motions. This assumption has been made to simplify the
notation.  Similar results can be obtained for
Markov processes and in this remark we indicate the essential arguments. \\
To be precise, assume that  the error processes $\{\boldsymbol{\varepsilon}(t) = (\varepsilon_1(t), \varepsilon_2(t))^\top| ~t\in[a, b]\}$   in model \eqref{moddisc} consist of two independent centered Gaussian processes with continuous covariance kernel given by
\begin{equation}\label{tri_kernel}
    K(s, t) = \mathbb{E}[\varepsilon_i(s) \varepsilon_i(t) ]
  =  u(s)v(t) \min\{q(s), q(t)\} \quad s, t \in [a, b] \, ,
\end{equation}
where $u(\cdot)$ and $v(\cdot)$ are  functions defined on the interval $[a, b]$ such that the function
$q(\cdot) = u(\cdot)/v(\cdot)$ is   positive and strictly increasing.  Kernels of the form \eqref{tri_kernel}
are called {\it triangular} kernels and a famous result in  \cite{doob1949heuristic} essentially shows that a Gaussian process is a Markov process if and only if its covariance kernel is triangular
\citep[see also][]{mehr1965certain}.
In this case model \eqref{moddisc} can be transformed into a model with an error process consisting of two independent Brownian motions
using the arguments given in Appendix $B$ of  \cite{detpepzhi2015}. More precisely,  define
\begin{align*}
 q(t) = \frac{u(t)}{v(t)}
\end{align*}
and consider the stochastic process
\begin{equation*}
\boldsymbol{\varepsilon} (t) = v(t) \boldsymbol
{ \tilde{\varepsilon}} ({q(t)}) ,
\end{equation*}
where $\{ \boldsymbol{\tilde{\varepsilon}} (\tilde{t}) = (\tilde{\varepsilon}_1({\tilde{t}})^\top , \tilde{\varepsilon}_2({\tilde{t}})) | \ \tilde{t} \in [\tilde{a},\tilde{b}] \}$  consists of two independent Brownian motions
on the interval $[\tilde{a},\tilde{b}]= [ q (a),   q (b)]$.  It now  follows
from  \cite{doob1949heuristic}
that  the process
$\{ \boldsymbol{ \varepsilon }(t)  = (\varepsilon_1(t), \varepsilon_2(t) )^\top  | \ t \in [a,b] \}$  consists of two independent centered Gaussian process on the interval $[a,b]$ with covariance kernel
\eqref{tri_kernel}. Consequently, if we consider the model
\be\label{mod1cont_transformed}
\boldsymbol{\tilde{Y}}(\tilde{t}) = \begin{pmatrix} \tilde{Y}_1(\tilde{t}) \\ \tilde{Y}_2(\tilde{t}) \end{pmatrix} = \mathbf{\tilde{F}}^\top (\tilde{t})\theta  + \mathbf{\Sigma}^{1/2}
 \boldsymbol{\tilde{\varepsilon}}(\tilde{t})
 \,  , \,  \,    \tilde{t} \in [q(a), q(b)],
\ee
and
$$
\tilde{\mathbf{F}}(\tilde{t}) =  \frac{\mathbf{F}(q^{-1}(\tilde{t}))}{v(q^{-1}(\tilde{t}))}~,~
 \boldsymbol{ \tilde{\varepsilon}} ({\tilde{t}})  = \frac{\boldsymbol{  \varepsilon} (q^{-1}(\tilde{t}))}{v(q^{-1}(\tilde{t}))} ~,~ \boldsymbol{ \tilde{Y}} ({\tilde t})  = \frac{\boldsymbol{ Y } (q^{-1}(\tilde{t}))}{v(q^{-1}(\tilde{t}))}\, ,
$$
the results obtained so far are applicable. Thus, a ''good'' estimator obtained for the parameter $\theta$ in model \eqref{mod1cont_transformed} is also a ''good estimator'' for the parameter $\theta$ in model \eqref{mod1cont} with error process consisting of two Gaussian processes with covariance kernel \eqref{tri_kernel}. Consequently,  we can derive the optimal estimator
for the parameter $\theta$ in the continuous time model \eqref{mod1cont}
with covariance kernel \eqref{tri_kernel} from the best linear unbiased estimator in the model given in \eqref{mod1cont_transformed} with Brownian motions by an application of Theorem \ref{thm1}. The resulting best linear unbiased estimator for $\theta$ in model \eqref{mod1cont} with triangular kernel \eqref{tri_kernel} is of the form
\begin{equation*}\label{asblue_transformed}
\hat{\theta}_{\rm BLUE}   = \mathbf{M}^{-1} \Big\{ \int_a^b \frac{\dot{\mathbf{F}}(t)v(t) - \mathbf{F}(t)v(t)}{\dot{u}(t)v(t)- u(t)\dot{v}(t)} \mathbf{\Sigma}^{-1}\,d\left(\frac{\mathbf{Y} (t)}{v(t)}\right) +   \frac{\mathbf{F}(a)\mathbf{\Sigma}^{-1} \mathbf{Y}(a)}{u(a)v(a)} \Big\},
\end{equation*}
where the minimum variance is given by
\begin{equation*}\label{cmatrix_inv_transformed}
\mathbf{M}^{-1} =\left(\int_a^b \frac{\left(\dot{\mathbf{F}}(t)v(t)- \mathbf{F}(t) \dot{v}(t)\right)\mathbf{\Sigma}^{-1}\left(\dot{\mathbf{F}}(t)v(t)- \mathbf{F}(t) \dot{v}(t)\right)^\top }{v^2(t) [\dot{u}(t)v(t)- u(t)\dot{v}(t)]} \,dt + \frac{ \mathbf{{F}}(a)\mathbf{\Sigma}^{-1} \mathbf{{F}}^\top (a)}{u(a)v(a)}\right)^{-1} \, .
\end{equation*}
}
 \end{remark}

\section{Optimal designs for comparing curves}
\label{sec4}
\def\theequation{4.\arabic{equation}}
\setcounter{equation}{0}

In this section we will derive optimal designs for comparing curves. The first part is devoted to a discretization of the
BLUE in the continuous time model. In the second part we develop
an optimality criterion to obtain efficient designs for the comparison of curves based on the discretized estimators.

\subsection{From the continuous to the discrete model }
\label{sec41}
To obtain a discrete  design
for $n$ observations at the points $a=t_1, \ldots , t_n$ from the continuous design derived in Section \ref{sec3},
we use a similar approach as in  \cite{detkonzhi2017} and  construct a discrete approximation of the stochastic integral in \eqref{asblue}. For this purpose we
 consider the linear estimator
\begin{eqnarray}\label{discrete_est}
\hat{\theta}_n &=& \mathbf{M}^{-1} \Big\{ \sum_{i=2}^n \mathbf{\Omega}_i \dot{\mathbf{F}}(t_{i-1})\mathbf{\Sigma}^{-1} (Y({t_i})-Y({t_{i-1}})) + \frac{\mathbf{F}(a)}{a}\mathbf{\Sigma}^{-1} Y_a \Big\} \\ \nonumber
&=& \mathbf{M}^{-1} \Big\{ \sum_{i=2}^n \mathbf{\Phi}_i  \mathbf{\Sigma}^{-1} (Y(t_i)-Y({t_{i-1}})) + \frac{\mathbf{F}(a)}{a}\mathbf{\Sigma}^{-1} Y_a \Big\},
\end{eqnarray}
were $a= t_1 <  t_2 < \ldots< t_{n-1}< t_n = b$,  $\mathbf{\Omega}_2, \ldots, \mathbf{\Omega}_n$ are $p \times p$ weight matrices and $\mathbf{\Phi}_2=\mathbf{\Omega}_2 \dot{\mathbf{F}}(t_{1}), \ldots, \mathbf{\Phi}_n = \mathbf{\Omega}_n \dot{\mathbf{F}}(t_{n-1})$ are $p\times 2$ matrices, which have to be chosen in a reasonable way. The matrix $\mathbf{M}^{-1}$ is given in \eqref{cmatrix_inv}.  To determine these weights in an ``optimal'' way we
first
derive a representation of the mean squared error between the best linear estimate \eqref{asblue}
in the continuous time model and
its discrete approximation \eqref{discrete_est}. The following result is a direct consequence of Ito's formula.

 \begin{lemma} \label{criterion-multi-dimension}
Consider the continuous time model \eqref{mod1cont}. If the assumptions of Theorem \ref{thm1} are satisfied, we have
\begin{align}
&\mathbb{E}_\theta[(\hat{\theta}_{\rm BLUE}  - \hat{\theta}_n)(\hat{\theta}_{\rm BLUE}  - \hat{\theta}_n)^\top ] = \mathbf{M}^{-1} \Big\{ \sum_{i=2}^n \int_{t_{i-1}}^{t_i} \big[ \dot{\mathbf{F}}(s) -  \mathbf{\Phi}_i \big] \mathbf{\Sigma}^{-1}\big[ \dot{\mathbf{F}}(s) - \mathbf{\Phi}_i   \big]^\top  \,ds \nonumber \\
&+ \sum_{i,j=2}^n \int_{t_{i-1}}^{t_i} \big[ \dot{\mathbf{F}}(s) - \mathbf{\Phi}_i \big]\mathbf{\Sigma}^{-1}\dot{\mathbf{F}}^\top (s) \,ds \, \theta \, \theta^\top
 \int_{t_{j-1}}^{t_j}  \dot{\mathbf{F}}(s)\mathbf{\Sigma}^{-1}\big[ \dot{\mathbf{F}}(s) - \mathbf{\Phi}_i \big]^\top  \,ds \Big\} \mathbf{M}^{-1}.
 \label{criterion-multi}
\end{align}
\end{lemma}

\bigskip

In the following we choose optimal $p\times 2$ matrices $\Phi_i=\mathbf{\Omega}_i \dot{\mathbf{F}}(t_{i-1})$ and design points $t_2, \ldots, t_{n-1}$ $(t_1=a, t_n=b)$, such that the linear estimate \eqref{discrete_est} is unbiased and    the mean squared error matrix in \eqref{criterion-multi} ``becomes small''.
An  alternative criterion is to replace the mean squared error $\mathbb{E}_\theta[(\hat \theta_{\rm BLUE}  - \hat \theta_n)(\hat \theta_{\rm BLUE}  - \hat \theta_n)^\top ]$  by the mean squared error
$$\mathbb{E}_\theta [(\hat{\theta}_n - \theta)(\hat{\theta}_n - \theta)^\top ]$$  between the   estimate $\hat \theta_n$ defined in \eqref{discrete_est} and the ``true'' vector of parameters.
The following result shows that
in the class of unbiased estimators both optimization problems  yield the same solution. The proof is similar to the proof of Theorem 3.1 in \citet{detkonzhi2017}.

\begin{theorem}\label{equivalence}
The estimator $\hat{\theta}_n$ defined in \eqref{discrete_est} is unbiased if and only if the identity
\begin{equation} \label{unbiasmul}
\mathbf{M}_0 = \int_a^b \dot{\mathbf{F}}(s)\mathbf{\Sigma}^{-1} \dot{\mathbf{F}}^\top (s) \,ds = \sum_{i=2}^n \Phi_i\mathbf{\Sigma}^{-1} \int_{t_{i-1}}^{t_i} \dot{\mathbf{F}}^\top (s) \,ds = \sum^n_{i=2}   \mathbf{\Phi}_i\mathbf{\Sigma}^{-1}(\mathbf{F}(t_i)-\mathbf{F}(t_{i-1}))^\top  ,
\end{equation}
is satisfied.
Moreover, for any linear unbiased estimator of the form $\tilde{\theta}_n = \int_a^b \mathbf{G}(s) dY_s $ we have
\begin{equation*} 
\mathbb{E}_\theta [(\tilde{\theta}_n - \theta)(\tilde{\theta}_n - \theta)^\top ] = \mathbb{E}_\theta [(\tilde{\theta}_n - \hat{\theta}_{\rm BLUE} )(\tilde{\theta}_n - \hat{\theta}_{\rm BLUE} )^\top ] + \mathbf{M}^{-1}.
\end{equation*}
\end{theorem}

\bigskip

In order to describe a solution in terms of optimal ``weights'' $\mathbf{\Phi}^*_i$ and design points $t^*_i$ we recall that the condition of unbiasedness of the estimate $\hat \theta_n$ in \eqref{discrete_est} is given by \eqref{unbiasmul} and
introduce the notation
\begin{align}
& \mathbf{B}_i = [\mathbf{F}(t_i) - \mathbf{F}(t_{i-1})]\mathbf{\Sigma}^{-1/2} / \sqrt{t_i-t_{i-1}},  \label{h1} \\
& \mathbf{A}_i = \mathbf{\Phi}_i\mathbf{\Sigma}^{-1/2} \sqrt{t_i-t_{i-1}}.\nonumber
\end{align}
It   follows from Lemma \ref{criterion-multi-dimension}
and Theorem  \ref{equivalence} that for an unbiased estimator $\hat \theta_n$ of the form \eqref{discrete_est}
the mean squared error has the representation
\be \label{crit}
 \mathbb{E}_\theta \big  [( \hat \theta_{\rm BLUE}  - \hat \theta_n)^\top  (\hat \theta_{\rm BLUE}  - \hat \theta_n)\big ]
  =
 -\mathbf{M}^{-1} \mathbf{M}_0 \mathbf{M}^{-1}   + \sum_{i=2}^n\mathbf{M}^{-1}  \mathbf{A}_i \mathbf{A}_i{^\top } \mathbf{M}^{-1}   ,
\ee
which has to be ``minimized'' subject to the constraint
\be \label{constraint}
\mathbf{M}_0 = \int^b_a \dot{\mathbf{F}} (s) \mathbf{\Sigma}^{-1} \dot{\mathbf{F}}^\top (s)ds=  \sum_{i=2}^n \mathbf{A}_i \mathbf{B}_i^\top.
\ee
The following result shows that a minimization with respect to  the weights $\mathbf{\Phi}_i$ (or equivalently $\mathbf{A}_i$) can actually be carried out with respect to
the Loewner ordering.

\begin{theorem} \label{thm3}
Assume that the assumptions of Theorem \ref{thm1} are satisfied and that the matrix
\begin{equation}\label{Bmatrix}
\mathbf{B} = \sum_{i=2}^{n} \mathbf{B}_i \mathbf{B}^\top_i =
\sum^n_{i=2} \frac {[\mathbf{F}(t_i) - \mathbf{F}(t_{i-1})]\mathbf{\Sigma}^{-1}[\mathbf{F}(t_i) - \mathbf{F}(t_{i-1})]^\top  }{t_i - t_{i-1}} ,
\end{equation}
is non-singular. Let  $\mathbf{\Phi}^*_2, \ldots, \mathbf{\Phi}^*_n$ denote $(p \times 2)$-matrices satisfying the equations
\be \label{eq1}
\mathbf{\Phi}^*_i = \mathbf{M}_0 \mathbf{B}^{-1} \frac {\mathbf{F}(t_i) - \mathbf{F}(t_{i-1})}{t_i - t_{i-1}}\qquad i=2,\ldots,n,
\ee
then $\mathbf{\Phi}^*_2, \ldots, \mathbf{\Phi}^*_n$ are optimal weight matrices minimizing $ \mathbb{E}_\theta[(\hat{\theta}_{\rm BLUE}
 - \hat{\theta}_n)(\hat{\theta}_{\rm BLUE} - \hat{\theta}_n)^\top ] $
 with respect to the Loewner ordering
among all unbiased estimators of the form \eqref{discrete_est}. Moreover, the variance of the resulting estimator $\hat\theta^*_n$ is given by
\begin{equation*}
\mathrm{Cov}(\hat\theta^*_n) = \mathbf{M}^{-1}\left\{\mathbf{M}_0 \mathbf{B}^{-1}\mathbf{M}_0 + \frac{1}{a} \mathbf{F}(a)\mathbf{\Sigma}^{-1}\mathbf{F}^\top(a)\right\}\mathbf{M}^{-1}
\end{equation*}
\end{theorem}

\textbf{Proof.} 
Let $v$ denote a $p$-dimensional vector and consider the problem of  minimizing
 the criterion
\begin{equation}\label{crit0}
v^\top \mathbb{E}_\theta[(\hat{\theta}_{\rm BLUE}  - \hat{\theta}_n)(\hat{\theta}_{\rm BLUE}
 - \hat{\theta}_n)^\top ]  v
 \end{equation}
  subject to the constraint \eqref{constraint}. Observing \eqref{crit} this yields the Lagrange function
\begin{equation}\label{lag_fun}
G_{v}(\mathbf{A}_1\ldots, \mathbf{A}_n) = - v^\top \mathbf{M}^{-1} \mathbf{M}_0 \mathbf{M}^{-1} v+ \sum^n_{i=2} (v^\top \mathbf{M}^{-1} \mathbf{A}_i\mathbf{A}_i^\top  \mathbf{M}^{-1}v) - \mbox{tr}\big\{\mathbf{\Lambda}(\mathbf{M}_0 - \sum_{i=2}^{n} \mathbf{A}_i \mathbf{B}_i^\top) \big\} ,
\end{equation}
where $\mathbf{A}_2, \ldots, \mathbf{A}_n$ are $(p\times 2)$-matrices and $\mathbf{\Lambda} = (\lambda_{k,\ell})^p_{k,\ell =1}$ is a $(p\times p)$-matrix of Lagrange multipliers. The
function $G_v$ is convex with respect to $\mathbf{A}_2,\ldots,\mathbf{A}_n$. Therefore, taking derivatives with respect to $\mathbf{A}_j$ yields as necessary and  sufficient  for the extremum  (here we use matrix differential calculus)
$$
2(\mathbf{M}^{-1}v)^\top\mathbf{A}_i\otimes (\mathbf{M}^{-1}v)^\top + \mbox{vec}\left\{\mathbf{\Lambda}\mathbf{B}_i\right\} = 0^\top_{2p} , \qquad i=2, \ldots, n  \, .
$$
Rewriting this system of linear equations in a $(p\times 2)$-matrix form gives
$$2\mathbf{M}^{-1}vv^\top \mathbf{M}^{-1} \mathbf{A}_i = - \mathbf{\Lambda} \mathbf{B}_i \qquad i =2, \ldots, n \, .$$
Substituting the expression in \eqref{constraint} and using the non-singularity of  the matrices $\mathbf{M}$ and $\mathbf{B}$ yields for the matrix of Lagrangian multipliers
$$\mathbf{\Lambda} = -2\mathbf{M}^{-1}vv^\top \mathbf{M}^{-1}\mathbf{M}_0\mathbf{B}^{-1}~, $$
which gives
\begin{equation}\label{iterim_result}
2\mathbf{M}^{-1}vv^\top \mathbf{M}^{-1} \mathbf{A}_i =  2\mathbf{M}^{-1}vv^\top \mathbf{M}^{-1}\mathbf{M}_0\mathbf{B}^{-1} \mathbf{B}_i \qquad i=2, \ldots , n \, .
\end{equation}
Note that one solution of \eqref{iterim_result} is given by
\begin{equation*}
\mathbf{A}^*_i = \mathbf{M}_0 \mathbf{B}^{-1}\mathbf{B}_i , \qquad i=2, \ldots, n
\end{equation*}
which does not depend on the vectors $v$. Therefore, the tuple of matrices $(\mathbf{A}^*_2, \ldots, \mathbf{A}^*_n)$ minimizes the convex
function $G_v$ in \eqref{lag_fun} for all $v\in \mathbb{R}^{p}$.\\
Observing the notations in \eqref{h1} shows that the optimal matrix weights are given by \eqref{eq1}. Moreover, these weights in \eqref{eq1}  do not depend on the vector $v$ either and
 provide a simultaneous
 minimizer of the criterion defined in \eqref{crit0} for all $v\in \mathbb{R}^p$. Consequently, the weights defined in \eqref{eq1} minimize $ \mathbb{E}_\theta[(\hat{\theta}_{\rm BLUE}
 - \hat{\theta}_n)(\hat{\theta}_{\rm BLUE} - \hat{\theta}_n)^\top ] $ under the unbiasedness constraint \eqref{constraint} with respect to the Loewner ordering. \hfill $\Box$

\bigskip

\begin{remark} \label{reminverse}
{\rm
If the matrix $\mathbf{B}$ in Theorem \ref{thm3} is singular, the optimal weights are not uniquely determined and we propose to replace the inverse $B$ by its Moore-Penrose inverse.}
\end{remark}

Note that for fixed design points $t_1, \ldots , t_n $ Theorem \ref{thm3} yields
universally optimal weights $\mathbf{\Phi}^*_2,\ldots,\mathbf{\Phi}^*_n$ (with respect to the Loewner ordering)
for  estimators of the form \eqref{discrete_est} satisfying \eqref{unbiasmul}. On the other hand, a further optimization with respect to the Loewner ordering with respect to the choice of the points $t_2,\ldots,t_{n-1}$ $(t_1=a, t_n=b)$ is not possible, and we have to apply a real valued optimality criterion for this purpose. In the following section, we will derive such a criterion which explicitly addresses the comparison of the regression curves from the
two groups introduced in Section \ref{sec2}.

\subsection{Confidence bands}
\label{sec42}

We return to the practical scenario of the two groups introduced in \eqref{indivmod}, where we now focus on the comparison of these groups on the interval $[a, b]$. \\
More precisely, consider the model introduced in \eqref{moddisc} and let $\hat{\theta}^*_n$ be the estimator   \eqref{discrete_est} with   optimal weights defined by \eqref{eq1} 
from  $n$ observations taken   at the time points $a=t_1<t_2 <\ldots < t_{n-1} <t_n = b$. Then this estimator is normally distributed with mean $\mathbb{E}[\hat{\theta}^*_n]=\theta$ and covariance matrix
$$\mbox{Cov}(\hat{\theta}^*_n)= \mathbf{M}^{-1}\left\{\mathbf{M}_0 \mathbf{B}^{-1}\mathbf{M}_0 + \frac{1}{a} \mathbf{F}(a)\mathbf{\Sigma}^{-1}\mathbf{F}^\top(a)\right\}\mathbf{M}^{-1}$$
where the matrices $\mathbf{M}, \mathbf{M}_0$ and $\mathbf{B}$ are given by \eqref{cmatrix}, \eqref{constraint}  and \eqref{Bmatrix},  respectively. Note that the covariance matrix depends on the time points $t_1, \ldots, t_n$ through the matrix $\mathbf{B}^{-1}$.
Moreover, using the estimator $\hat{\theta}^*_n$ the prediction of the difference of a fixed time point $t\in [a,b]$ satisfies
\begin{equation*}
(1, -1) \mathbf{F}^\top(t) \hat{\theta}^*_n - (1, -1) \mathbf{F}^\top(t) \theta \sim \mathcal{N}_p(0, h(t; t_1, \ldots, t_n))~,
\end{equation*}
where
\begin{equation*}
    h(t; t_1, \ldots, t_n) = (1, -1)\mathbf{F}^\top \mathbf{M}^{-1}\left\{\mathbf{M}_0 \mathbf{B}^{-1}\mathbf{M}_0 + \frac{1}{a} \mathbf{F}(a)\mathbf{\Sigma}^{-1}\mathbf{F}^\top(a)\right\}\mathbf{M}^{-1}\mathbf{F}(t) (1, -1)^T \, .
\end{equation*}
We now use this result and the results of \cite{gsteiger2011} to obtain a simultaneous confidence band for the difference of the two curves. More precisely, if the interval $[a,b]$ is the range where the two curves should be compared, the simultaneous confidence band is defined as follows.  Consider the statistic
\begin{equation*}\label{confband}
 \hat T = \sup_{t \in [a,b]} \ \frac {|(1, -1) \mathbf{F}^\top(t) \hat{\theta}^*_n - (1, -1) \mathbf{F}^\top(t) \theta |}{  \{h(t; t_1, \ldots, t_n)\}^{1/2}} ,
\end{equation*}
and define $D$ as the $(1-\alpha)$-quantile of the corresponding distribution, that
is
\begin{equation*}
\mathbb{P}(\hat T \leq D) = 1-\alpha.
\end{equation*}
Note that \citet{gsteiger2011} propose the parametric bootstrap for choosing the critical value $D$.
Define
\begin{eqnarray*}
u (t; t_1, \ldots, t_n) &= &
(1, -1) \mathbf{F}^\top(t) \hat{\theta}^*_n
+ D \cdot {  \{h(t; t_1, \ldots, t_n)\}^{1/2}} ,\\
l (t; t_1, \ldots, t_n) &= &
(1, -1) \mathbf{F}^\top(t) \hat{\theta}^*_n
- D \cdot {  \{h(t; t_1, \ldots, t_n)\}^{1/2}},
\end{eqnarray*}
then the confidence band  for the difference of the two regression functions is defined by
\begin{eqnarray}\label{conf}
{\cal C}_{1-\alpha}
= \big \{ g: [a,b] \to \mathbb{R}~|~
l (t; t_1, \ldots, t_n) \leq g(t ) \leq  u (t; t_1, \ldots, t_n) \mbox{ for all } t \in [a,b] \big \}~.
\end{eqnarray}
Consequently, good time points $t_1=a< t_2< \ldots<t_{n-1}, t_n=b$ should minimize  the width
$$
u (t; t_1, \ldots, t_n) - l (t; t_1, \ldots, t_n) =
2  \cdot D \cdot {  \{h(t; t_1, \ldots, t_n)\}^{1/2}}~
$$
of this band  at each $t \in [a,b ]$.
 As  
 this is only possible in rare circumstances,  we propose to minimize an $L_p$-norm of the function $h(\cdot ; t_1\ldots, t_n)$ as a design criterion, that is
\begin{equation}\label{lp}
\Phi_p(t_1, \ldots, t_n) = \| h(\cdot; t_1\ldots, t_n) \|_p := \Big( \int_a^b [h(t; t_1\ldots, t_n)]^p \Big)^{1/p} \,dt, \quad 1 \leq p \leq \infty,
\end{equation}
where the case $p=\infty$ corresponds to the maximal deviation  $$
\| h(\cdot; t_1\ldots, t_n) \|_{\infty} = \sup_{t \in [a,b]} |h(t; t_1\ldots, t_n)|.
$$
Finally, the optimal points $a=t^*_1 < t_2^* < \ldots < t_n^*=b$ (minimizing \eqref{lp}) and the corresponding weights derived in Theorem \ref{thm3} provide the optimal linear unbiased estimator of the form \eqref{discrete_est} (with the corresponding optimal design).

\begin{example}
{\rm
We now conclude this section by considering the cases of  no common and common parameters, respectively.\\
(a) If we are in the situation of Example \ref{ex1}
(no common parameters), the regression function $\mathbf{F}^\top(t)$ is of the form in \eqref{modcont_sep1} and
the variance of the prediction of the difference at a fixed point $t\in[a,b]$ reduces to
$$h(t; t_1, \ldots, t_n) =  (f^\top_1(t), - f^\top_2(t)) \mathbf{M}^{-1}\left\{\mathbf{M}_0 \mathbf{B}^{-1}\mathbf{M}_0 + \frac{1}{a} \mathbf{F}(a)\mathbf{\Sigma}^{-1}\mathbf{F}^\top(a)\right\}\mathbf{M}^{-1} (f^\top_1(t), - f^\top_2(t))^\top. $$
The corresponding design criterion is given by
\begin{equation}\label{phi_p}
\Phi_p \big(t_1, \ldots, t_n\big) = \|(f^\top_1, - f^\top_2) \mathbf{M}^{-1}\left\{\mathbf{M}_0 \mathbf{B}^{-1}\mathbf{M}_0 + \frac{1}{a} \mathbf{F}(a)\mathbf{\Sigma}^{-1}\mathbf{F}^\top(a)\right\}\mathbf{M}^{-1}(f^\top_1, - f^\top_2)^\top\|_p \, .
\end{equation}
(b) If we are in the situation of  Example \ref{ex2}
(common parameters), the regression function $\mathbf{F}^\top(t)$ is given by \eqref{commonF} and the variance of the prediction of the difference at a fixed point $t\in[a, b]$ reduces to
$$h(t; t_1, \ldots, t_n) =  (0, \tilde{f}^\top_1(t), - \tilde{f}^\top_2(t)) \mathbf{M}^{-1}\left\{\mathbf{M}_0 \mathbf{B}^{-1}\mathbf{M}_0 + \frac{1}{a} \mathbf{F}(a)\mathbf{\Sigma}^{-1}\mathbf{F}^\top(a)\right\}\mathbf{M}^{-1}  (0, \tilde{f}^\top_1(t), - \tilde{f}^\top_2(t)) ^\top \, .  $$
The corresponding design criterion is given by
\begin{equation*}
\Phi_p \big(t_1, \ldots, t_n\big) = \|(0^\top_{p_0}, \tilde{f}^\top_1, - \tilde{f}^\top_2) \mathbf{M}^{-1}\left\{\mathbf{M}_0 \mathbf{B}^{-1}\mathbf{M}_0 + \frac{1}{a} \mathbf{F}(a)\mathbf{\Sigma}^{-1}\mathbf{F}^\top(a)\right\}\mathbf{M}^{-1}(0^\top_{p_0}, \tilde{f}^\top_1, - \tilde{f}^\top_2)^\top\|_p \, .
\end{equation*}
}
\end{example}

\section{Numerical Examples} \label{sec5}
In this section the methodology is illustrated  in examples by means of a simulation study. To be precise,  we consider the  regression model \eqref{moddisc}, where the matrix
$\mathbf{F} (t)$
is given by \eqref{modcont_sep1} corresponding to the case that
the regression function do not share common parameters, see  Section \ref{sec31} for more details.
In this case the corresponding bounds for the confidence band
is  given by \eqref{conf}, where
\begin{eqnarray*}
u (t; t_1, \ldots, t_n) &= & (\hat{\theta}^{*(1)}_n)^\top f^{(1)} (t) -  (\hat{\theta}^{*(2)}_n)^\top f^{(2)} (t)
+ D \cdot {  \{h(t; t_1, \ldots, t_n)\}^{1/2}} ,\\
l (t; t_1, \ldots, t_n) &= & (\hat{\theta}^{*(1)}_n)^\top f^{(1)} (t) -  (\hat{\theta}^{*(2)}_n)^\top f^{(2)} (t)
- D \cdot {  \{h(t; t_1, \ldots, t_n)\}^{1/2}},
\end{eqnarray*}
and
$\hat{\theta}^*_n
= ( (\hat{\theta}^{*(1)}_n)^\top, (\hat{\theta}^{*(2)}_n)^\top)^\top $ is the estimator \eqref{discrete_est} with  optimal weights defined in \eqref{eq1}.
The design space  is given by the interval
$[a, b] = [1, 10]$,
and  we consider
three choices for the functions
 $f_1$ and $f_2$ in the matrix
 \eqref{modcont_sep1}, that is
\begin{eqnarray}
\nonumber
    f_A(t) &=& (t, \sin(t), \cos(t))^\top  \, , \\
    \label{fABC}
     f_B(t) &=& (t^2, \cos(t), \cos(2t))^\top  \, , \\
  f_C(t) &=& \big(t, \log(t), \frac{1}{t}\big)^\top  \,  . \nonumber
\end{eqnarray}
To model  the dependence between the two groups we use the covariance matrix
$$\mathbf{\Sigma} = \begin{pmatrix}1 &  \varrho \\ \varrho & 1 \end{pmatrix} \, ,
$$
in \eqref{moddisc},
where the correlations
are chosen as $\varrho=  0.2, 0.5,  0.7$.
Following the discussion in Section \ref{sec41} we focus on the comparison of the
regression curves for the two groups and derive optimal designs, minimizing the criterion  $ \Phi_\infty$ defined in
\eqref{phi_p}. As result, we obtain simultaneous
 confidence bands  with a smaller maximal width for the difference of the curves describing the relation in the two groups. We can obtain similar results
different  values $p \in (0,\infty) $ in \eqref{phi_p} but for the sake of brevity we concentrate on the  criterion  $\Phi_{\infty}$ which is
probably  also the easiest to interpret  for practitioners.

We denote by $\hat{\theta}_{n}^*$  the linear unbiased estimator derived in Section \ref{sec4}.  For  each of the combinations of regression functions containing two different functions defined in \eqref{fABC}, the
optimal weights have been found by  Theorem \ref{thm3} 
and the optimal design points  $t_{i}^*$
are determined minimising the criterion $\Phi_\infty$ defined in \eqref{phi_p}.
For the numerical  optimisation
 the Particle Swarm Optimisation (PSO) algorithm is used  \citep[see, for example,][]{Clerc2006}
 assuming a sample size of four observations
in each group, that is, $n=4$.  Furthermore, the uniform design used in the following calculations  is the design which has four equally spaced design points in the interval $[1, 10]$. The $\Phi_{\infty}$-optimal design points
minimizing the criterion
criterion in \eqref{phi_p} are given in Table \ref{tab1} for all  combinations of models and  correlations
under consideration. Note  that for each model  the corresponding optimal design points change for different values of correlation $\varrho$.

\begin{table}[h]
\caption{\label{tab1} \emph{Optimal designs points on the interval $[1,10]$ for the  estimator $\hat{\theta}_{n}^*$  in \eqref{discrete_est} minimizing the criterion $\Phi_{\infty}$  in \eqref{phi_p}. Different
 correlations  $\varrho = 0.2, 0.5, 0.7$ and different regression functions defined in  \eqref{fABC} are considered.}}
\vspace{0.5cm}
\centering
\begin{tabular}{|c||c|c|c|}
\hline
 & \multicolumn{3}{|c|}{correlation} \\
\hline
models & $\varrho =0.20$
& $\varrho = 0.50$ &  $\varrho = 0.70$\\
\hline
$f_1 = f_A$ \& $f_2 = f_B$ &  [1, 1.59, 3.93, 10] & [1, 1.62, 3.91, 10] & [1, 1.74, 7.99, 10]\\
\hline
 $f_1=f_A$ \& $f_2= f_C$  &  [1, 3.46, 9.60, 10] & [1, 2.86, 8.83, 10] & [1, 2.61, 3.52, 10]\\
\hline
$f_1 = f_B$ \& $f_2= f_C$  &[1, 2.20 , 6.25, 10]&[1, 1.62, 3.98 ,  10] & [1, 2.85, 6.29 ,  10] \\
\hline
\end{tabular}
\end{table}

\begin{figure}[!htbp!]
\centering\caption{{\it Confidence bands for the difference of the regression functions (solid grey line) on the basis of an optimal (solid lines) and uniform design (dashed lines). Left panel: $\varrho=0.2$. Middle panel: $\varrho = 0.5$. Right panel: $\varrho = 0.7$. First row: model with $f_1=f_A$ and $f_2=f_B$. Second row: model with $f_1=f_A$ and $f_2=f_C$. Third row: model with $f_1=f_B$ and $f_2=f_C$.}}
	\label{fig3_1}
  \begin{minipage}{.3\textwidth}
    \centering				\includegraphics[width=.8\textwidth, angle =270]{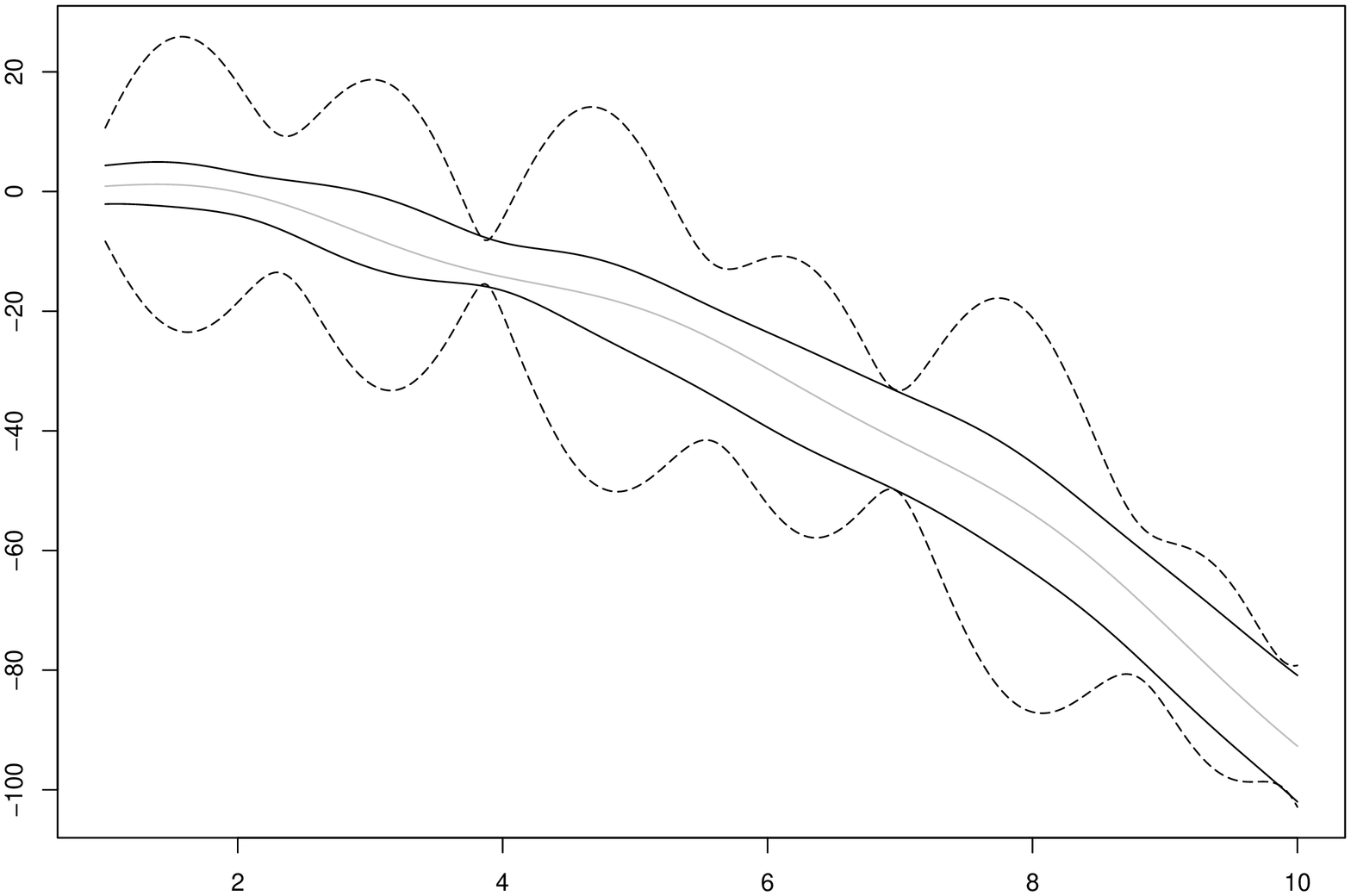} \\
		\includegraphics[width=.8\textwidth, angle =270]{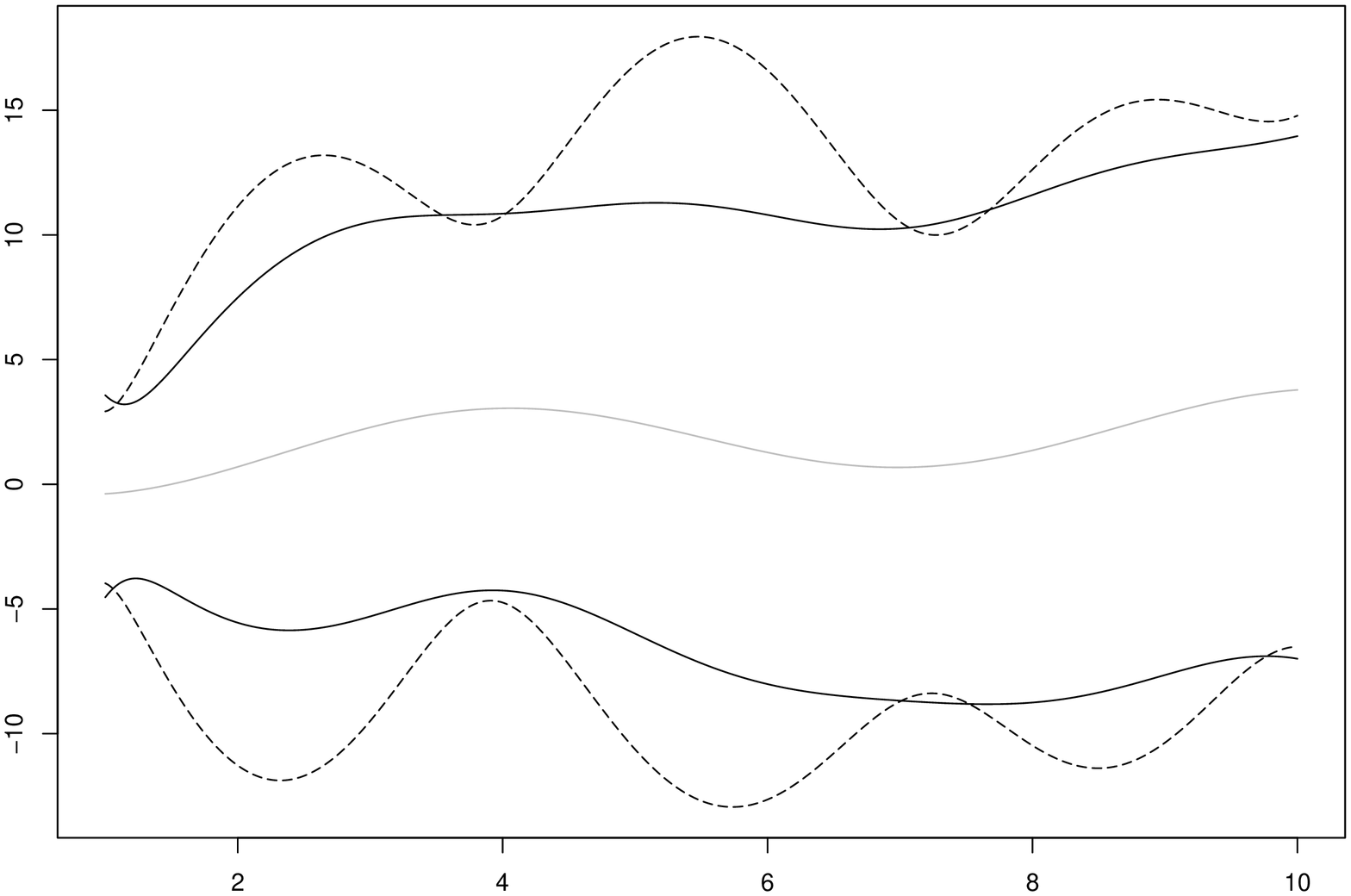} \\
		\includegraphics[width=.8\textwidth, angle =270]{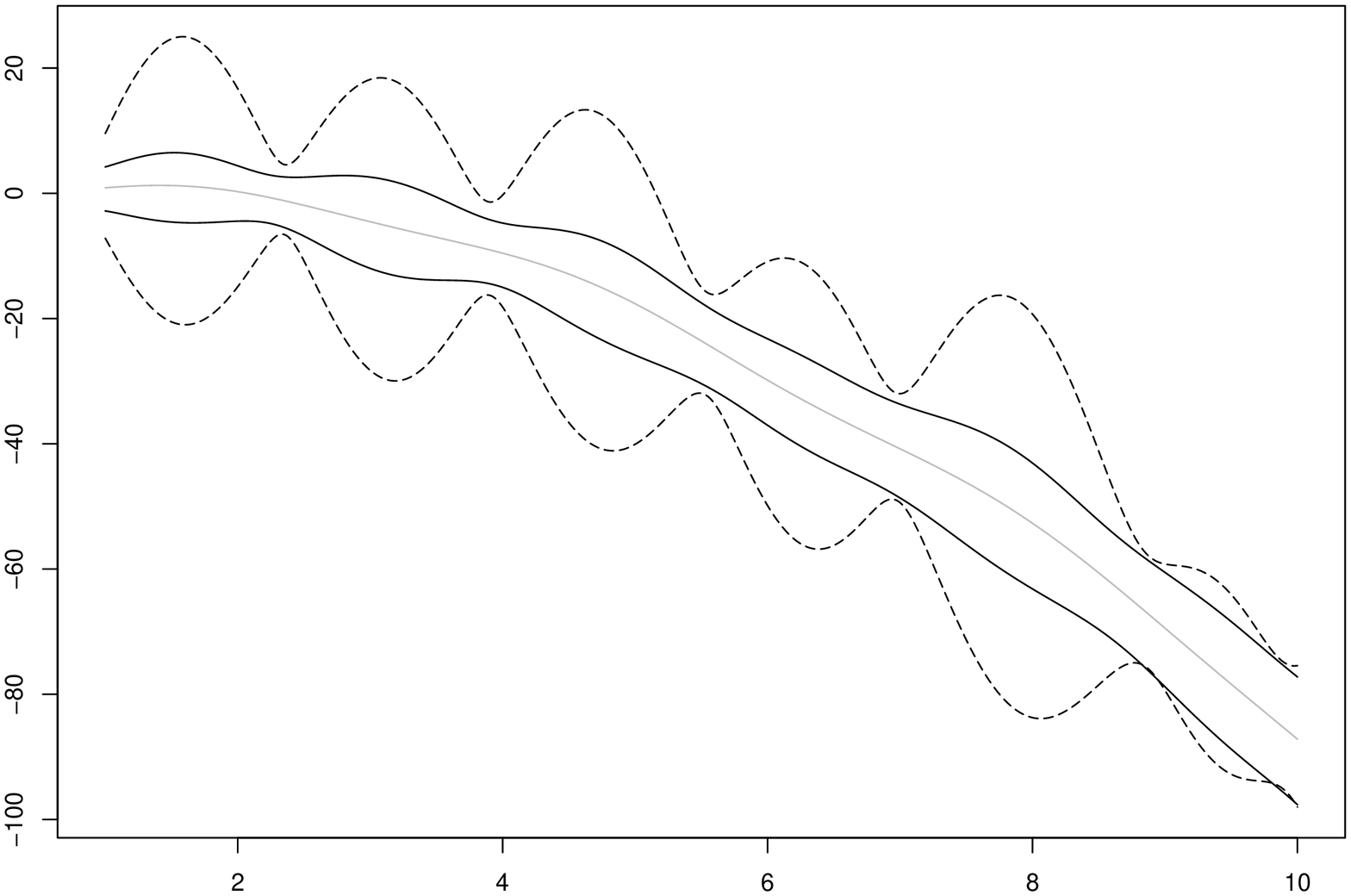} \\
		\end{minipage}%
		\begin{minipage}{.3\textwidth}
 \centering
		\includegraphics[width=.8\textwidth, angle =270]{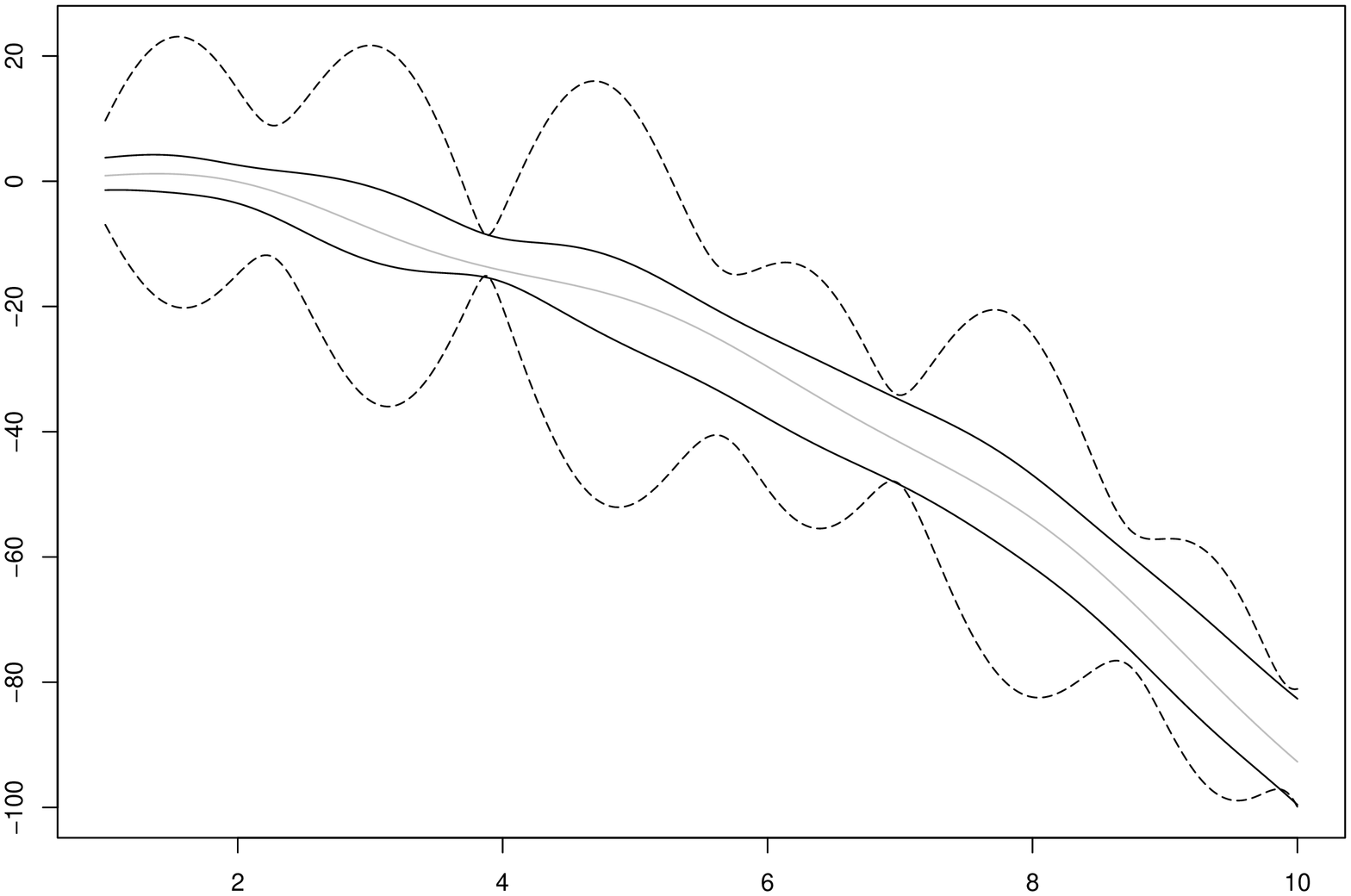} \\
		\includegraphics[width=.8\textwidth, angle =270]{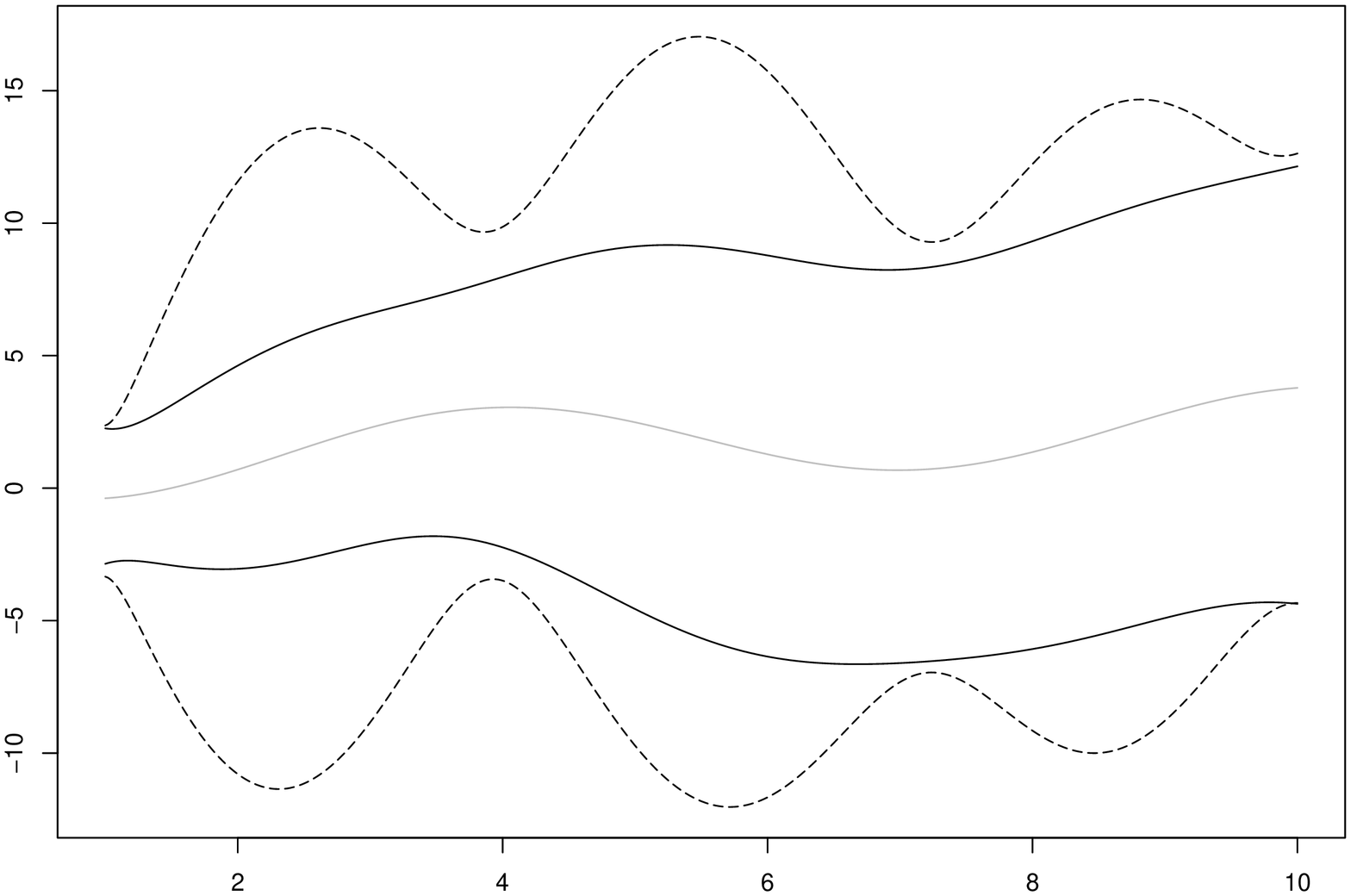} \\
		\includegraphics[width=.8\textwidth, angle =270]{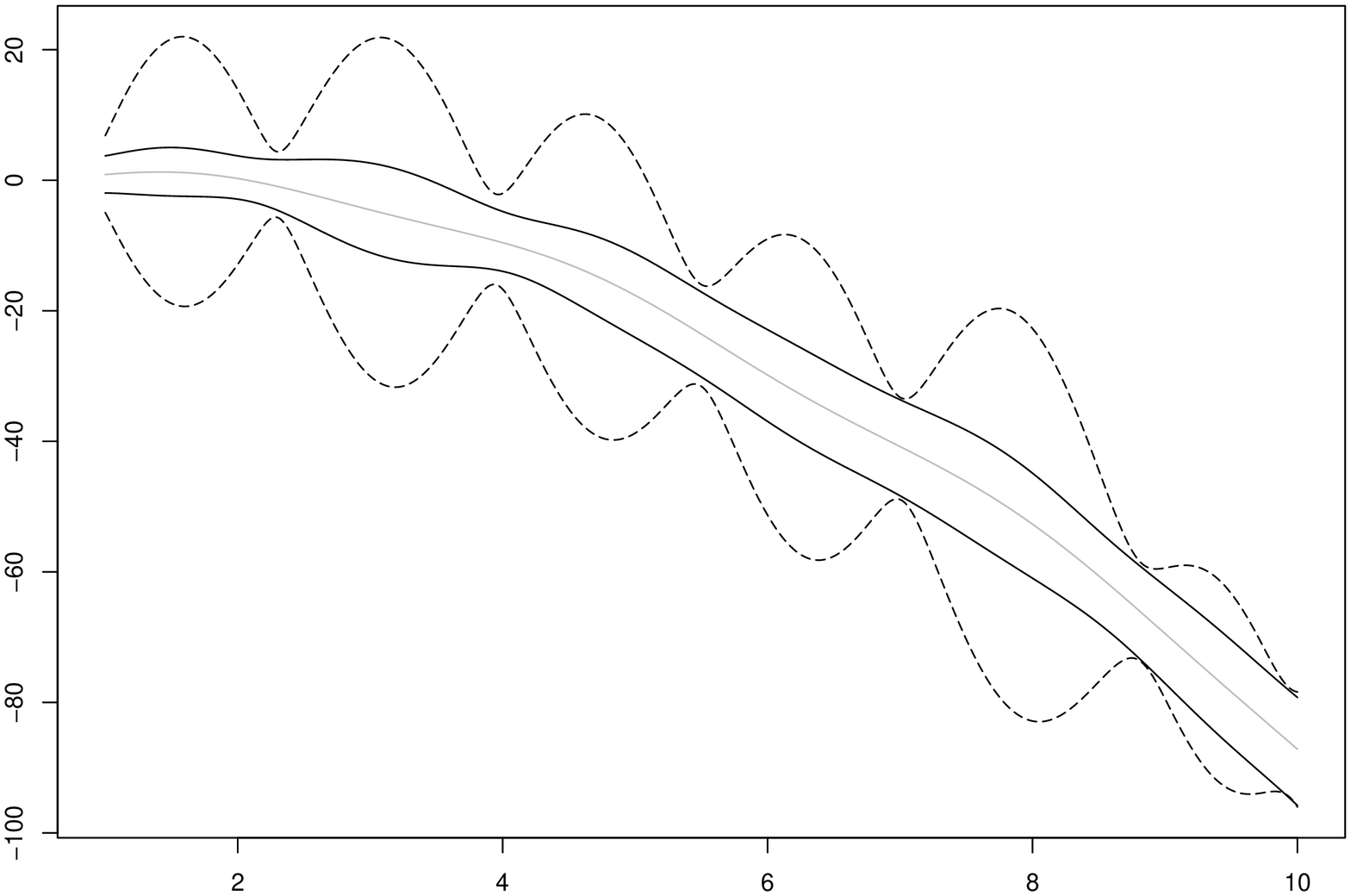}
		\end{minipage}
				\begin{minipage}{.3\textwidth}
 \centering
		\includegraphics[width=.8\textwidth, angle =270]{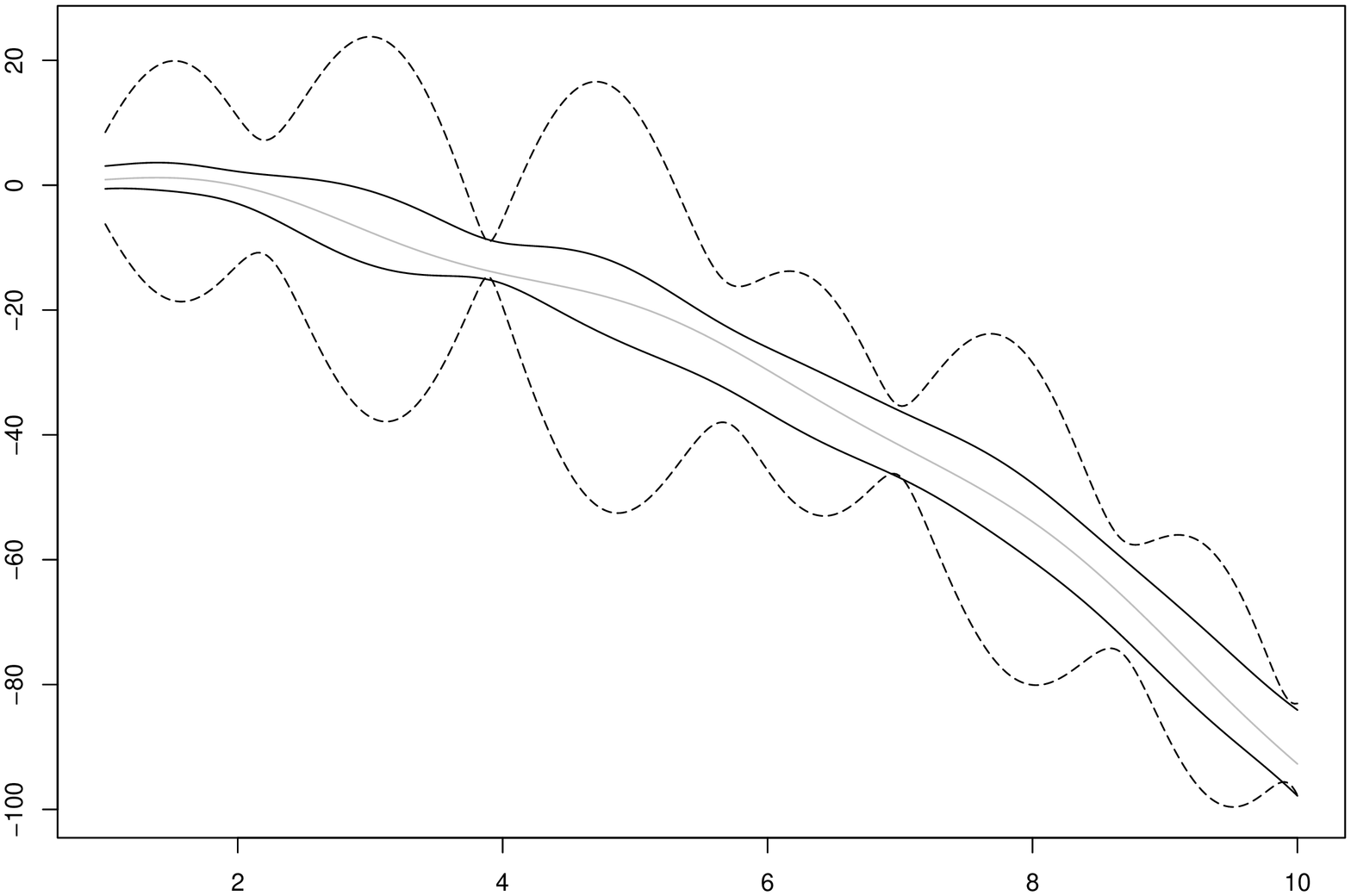} \\
		\includegraphics[width=.8\textwidth, angle =270]{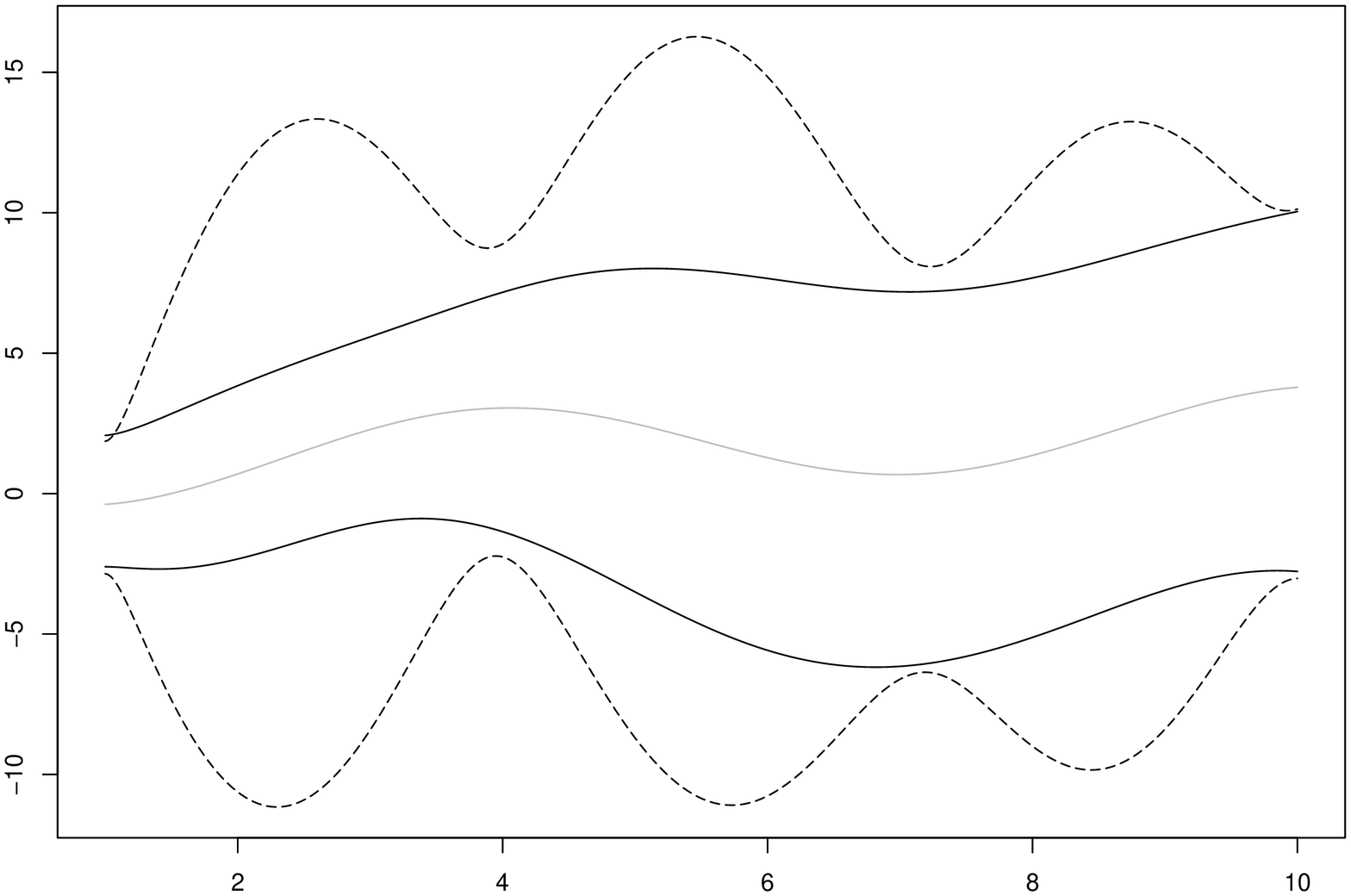} \\
		\includegraphics[width=.8\textwidth, angle =270]{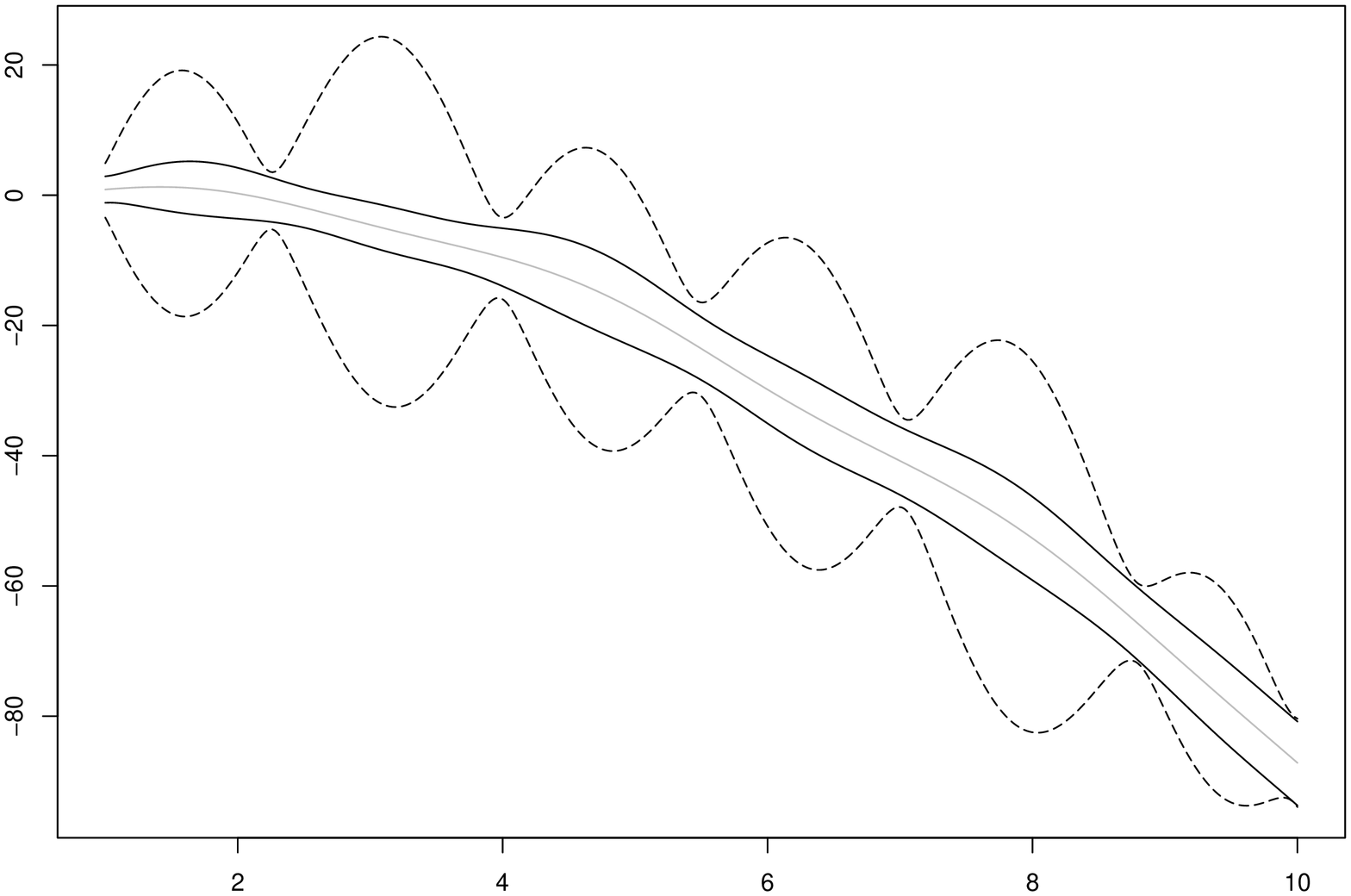}
		\end{minipage}
\end{figure}

In order to investigate the impact of the optimal design on the structure of the
confidence bands we have performed a small simulation study simulating confidence bands
for the difference of the regression functions. The vector of parameter values used for each model is  $\theta = ({\theta^{(1)}}^\top, {\theta^{(2)}}^\top)^\top = (1, 1, 1, 1, 1, 1)^\top$. In Figure \ref{fig3_1} we display
the  averages of uniform confidence bands
defined in \eqref{conf}
under the uniform and optimal design
calculated by $100$ simulation runs.

The left, middle and right columns
show the results  for the correlation $\varrho = 0.2$, $\varrho=0.5$ and $\varrho = 0.7$, respectively, while the rows correposnd to different combinations for  the functions $f_1$ and $f_2$ (first row: $f_1=f_A$, $f_2=f_B$, middle row:  $f_1 = f_A$, $f_2= f_C$ and last row  $f_1= f_B$, $f_2=f_C$).
In each graph, the confidence bands from the $\Phi_{\infty}$-optimal or the uniform design are plotted separately using the solid and dashed lines respectively, along with the plot for the true difference $f_1^\top(t)\theta^{(1)} - f_2^\top(t)\theta^{(2)}$ (solid grey lines).

We  observe, that in all cases under considerations the use of  $\Phi_{\infty}$-optimal designs
yields a clearly visible improvement compared to the  uniform design.   The maximal width of the confidence band is reduced substantially. Moreover,  the bands from the $\Phi_{\infty}$-optimal designs are nearly uniformly more narrow than the bands based on the  uniform design. Even more importantly,
the confidence bands based on the $\Phi_{\infty}$-optimal design
show a similar structure as the true differences, while the confidence bands from  the uniform design oscillate.

A comparison of  the left, middle and  right columns in  Figure \ref{fig3_1}
shows that  the maximum width for the confidence bands based on  the optimal design
decreases with increasing (absolute) correlation $\varrho$. This effect is not visible for the confidence bands  based on the uniform design. For example, for the middle row of Figure \ref{fig3_1}, which corresponds to the case $f_1 = f_A$ and $f_2 = f_C$, the maximum width of the confidence bands based on the equally spaced design points even seem to increase.\\

Table \ref{tab2} presents the values of the criterion  $\Phi_{\infty}$ in \eqref{phi_p} for the different scenarios  and confirms  the conclusions drawn from the visual inspection of the confidence bands plots.
We observe that the use of the optimal design points reduces the maximum width of the confidence bands substantially. Moreover, for the optimal design the maximum width becomes smaller with increasing (absolute) correlation. On the other hand this monotonicity cannot be observed in all cases for the uniform designs.

Summarizing, the use of the proposed $\Phi_{\infty}$-optimal design   improves statistically  inference   substantially reducing the maximum variance of the difference of the two estimated regression curves.   Moreover, simultaneous estimation in combination with a     $\Phi_\infty$-optimal design  yields a further reduction of the maximum width of the confidence bands, thus providing a more precise inference for the difference of the curves describing the relation between $t$ and the responses in the two groups.

\begin{table}[h]
\caption{\label{tab2} \emph{Values of the criterion  $\Phi_\infty$ for the optimal and uniform design with four observations in each group in the interval $[1,10]$. The error process is given by a two independent Brownian motions with correlation $\varrho = 0.2, 0.5, 0.7$ between the groups, respectively.}}
\vspace{0.5cm}
\centering
\begin{tabular}{|c|c|r|r|r|}
\hline
\multicolumn{2}{|c|}{} & \multicolumn{3}{|c|}{Correlation} \\
\hline
  Models & Design & $\varrho = 0.2$ & $\varrho = 0.5$ & $\varrho = 0.7$ \\
\hline
\multirow{2}{*}{$f_1= f_A$ \& $f_2 = f_B$}   & optimal &14.79 & 9.44 & 6.09 \\
\cline{2-5}
&  uniform  & 141.87 & 142.59 & 148.74    \\
\hline
\hline
\multirow{2}{*}{$f_1= f_A$ \& $f_2 = f_C$} & optimal & 16.00 & 10.00 & 6.60  \\
\cline{2-5}
  &uniform  & 33.32  & 29.10 & 25.66    \\
\hline
\hline
\multirow{2}{*}{$f_1= f_B$ \& $f_2 = f_C$} & optimal &14.71 & 9.53 &5.99 \\
\cline{2-5}
&  uniform  & 147.27  &127.19 &  115.07 \\
\hline
\end{tabular}
\end{table}


\bigskip
\bigskip

{\bf Acknowledgements}
This research was partially supported by the Collaborative Research Center `Statistical modeling
of nonlinear dynamic processes' ({\it Sonderforschungsbereich 823, Teilprojekt C2}).\\
On behalf of all authors, the corresponding author states that there is no conflict of interest. 
%

 \bibliography{opt_signed_designsa}

\end{document}